\documentclass[review]{elsarticle}

\def\bdk{\begin{description} \itemsep=-\parsep \itemindent=-0.9 cm}
\def\edk{\end{description}}
\def \leaderfill{\leaders\hbox to 0.35em{\hss$\cdot$\hss}\hfill}

\usepackage[fleqn]{amsmath}
\usepackage{mathrsfs}
\usepackage{bm}
\usepackage{amssymb}
\usepackage{color}
\usepackage{xcolor}
\usepackage{graphicx}
\usepackage{amsthm}
\usepackage{cases}
\usepackage{caption,array,threeparttable}
\usepackage{subcaption}
\usepackage{geometry}
\usepackage{epstopdf}
\usepackage[figuresright]{rotating}
\usepackage{booktabs}
\usepackage{longtable}
\usepackage{wasysym}
\usepackage{arydshln}
\usepackage{multirow}
\usepackage{float}
\usepackage{nicematrix}
\usepackage[english]{babel}
\usepackage{enumerate}
\usepackage{extarrows}
\usepackage{algorithm}
\usepackage{supertabular}
\usepackage[titletoc]{appendix}
\usepackage{algpseudocode}
\usepackage{array}
\usepackage{tabularx}
\usepackage{pdflscape}

\newtheorem{theorem}{Theorem}[section]

\newtheorem{proposition}{Proposition}[section]

\usepackage[font=small,labelfont=bf,labelsep=none]{caption}
\usepackage{natbib}
\biboptions{authoryear,round}

\allowdisplaybreaks[4]
\journal{Journal of LaTeX Templates}

\usepackage{hyperref}

\begin{document}

\begin{frontmatter}

\title{A new matrix-variate integer-valued autoregressive process with matrical negative binomial thinning}

\author[b]{Chunhan Liu}
\author[a]{Nuo Xu}
\author[b]{Kai Yang\corref{mycorrespondingauthor}}
\cortext[mycorrespondingauthor]{Corresponding author}
\ead{yangkai@ccut.edu.cn}

\address[a]{School of Mathematics, Jilin University, China}
\address[b]{School of Mathematics and Statistics, Changchun University of Technology, Changchun 130012, China}

\begin{abstract}
To address the overdispersion problem in matrix-variate integer-valued time series data arising in sociology, medicine, and related fields, this paper proposes a matrix integer-valued autoregressive model based on the negative binomial thinning operator. 
By introducing left and right matricial negative binomial thinning operators, the proposed model not only preserves the integer-valued nature of the data but also effectively handles overdispersion. The probabilistic and statistical properties of the proposed model are systematically investigated. 
Two estimation methods, namely projection estimation and iterative conditional least squares estimation are developed, and the corresponding asymptotic theories are established. 
Simulation studies provide concrete numerical results to evaluate the finite-sample performance of the estimators. 
Real data analysis demonstrates that the proposed model outperforms both the continuous matrix autoregressive model and the multivariate integer-valued autoregressive model in fitting matrix-variate integer-valued time series, while also proving effective in accommodating overdispersed count data.
\end{abstract}


\end{frontmatter}

\section{Introduction}

With the rapid advancement of big data and artificial intelligence technologies, data collection and storage capabilities have significantly improved, and the structure of observed data has become increasingly complex. Traditional point-valued data are simple in structure and easy to handle, but often fail to fully capture the true state of complex systems. To overcome the informational limitations of one-dimensional point data, some scholars have extended data to vectors or high-dimensional data consisting of different indicators for the same objects (\cite{Lam2012,Chen2013,Wang2022}). Other researchers have developed cross-sectional data (\cite{Kim2021}) or panel data (\cite{Jiang2023}) obtained by observing the same indicator across different regions. Although these two types of data integrate multiple indicators or multiple regions, they still have clear application boundaries, as their research objects are always limited to a single entity. With the acceleration of globalization and the strengthening of regional linkages, research in fields such as economics, finance, criminology, medicine, and insurance increasingly requires breaking through the limitation of a single entity to achieve cross-dimensional observation of multiple entities and multiple indicators. Consequently, researchers have attempted to integrate these two types of data for joint observation, giving rise to matrix data.

For the macroeconomic monitoring of major global economies, to compare the economic performance of China, the United States, and the United Kingdom, researchers typically need to construct matrix data organized by countries and economic indicators. The rows of the matrix correspond to the three countries, and the columns correspond to three core indicators: GDP growth rate, CPI inflation rate, and IPI growth rate. The value in each cell accurately reflects the performance of a specific country on a specific indicator, facilitating both horizontal comparisons of the same indicator across different countries and vertical analysis of the correlations among multiple indicators within the same country.

Early research on matrix data primarily focused on the distribution theory of random matrices (\cite{Chikuse1981}). Research during this period mainly aimed to characterize the probability distributions of matrix data, providing a theoretical foundation for subsequent modeling, but did not yet address the modeling or practical application of matrix data. It was not until the 21st century that the research focus on matrix data gradually shifted from distribution theory to modeling methods. \cite{Li2010} proposed a dimension folding method for matrix data dimensionality reduction, which effectively solves classification problems for matrix data and overcomes the limitation of traditional dimensionality reduction methods that destroy data structure, providing a new approach for the dimensionality reduction and classification of matrix data. \cite{Kong2012} defined a sparse matrix graphical model based on the matrix normal distribution, effectively capturing row and column dependencies in matrices, achieving better interpretability than traditional sparse vector graphical models. \cite{Zhang2014} proposed a nuclear norm regularized matrix regression method, introducing a regularization term to address the overfitting problem in high-dimensional matrix regression. \cite{Ding2018} proposed regression models and estimation methods for matrix-response regression under different scenarios, accounting for complex relationships among variables. \cite{Kong2020} proposed a bilinear transformation matrix regression model for high-dimensional matrix response data, reducing model complexity and improving estimation efficiency. \cite{Chang2023} proposed a matrix time series modeling method based on tensor CP decomposition, effectively addressing the over-parameterization problem in modeling. \cite{Bettache2025} proposed low-rank prediction methods for two-sided matrix regression models using rank penalty and nuclear norm penalized least squares.

Research on matrix data exhibits the following characteristics: first, it emphasizes the preservation of data structure; Compared to traditional vectorization approaches, matrix data modeling can fully utilize the row and column structural information of the data.
Second, it incorporates techniques such as dimensionality reduction and regularization to address the curse of dimensionality caused by high-dimensional matrix data. Third, its application scenarios are continuously expanding, from initial theoretical research in statistics to various practical fields such as economics, meteorology, and computer vision. However, research on matrix data still has some shortcomings. On the one hand, modeling methods for high-dimensional matrix data still need improvement, especially when both the row and column dimensions are large, the computational complexity of models is high, and the efficiency of parameter estimation needs to be enhanced. On the other hand, the dependence structures between rows and columns of matrix-variate data cannot be replaced by traditional vector data.

With the continuous advancement of data collection technology, matrix-type data are frequently observed sequentially over time, forming matrix time series. Matrix time series possess both the structural characteristics of matrix data and the dynamic dependencies of time series, enabling them to describe the evolution patterns of multiple objects and multiple indicators over time. Compared to traditional vector time series, they are better suited for characterizing the dynamic evolution of complex systems. For instance, in finance, daily data of multiple financial indicators for multiple stocks, or monthly data of multiple macroeconomic indicators across different countries; in environmental science, hourly data of multiple pollutant concentrations at multiple monitoring stations are all typical examples of matrix time series. The concept of matrix time series was first introduced by Walden and Serroukh (\cite{Walden2002}), who applied wavelet analysis methods to matrix time series analysis and constructed various matrix-valued filter models, marking the formal beginning of this field. \cite{Samadi2014} proposed the first-order matrix autoregressive model MAR(1), which breaks through the limitations of the traditional vector autoregressive model and lays the foundation for the important branch of matrix autoregressive models.
\cite{Chen2021} proposed an innovative matrix autoregressive model with a bilinear structure, which not only preserves the original matrix structure of the data but also achieves effective dimensionality reduction through the bilinear structure. In the same year, \cite{Hsu2021} incorporated spatial dependencies and proposed a matrix spatiotemporal autoregressive model. \cite{Mircetic2022} applied the matrix autoregressive model to demand forecasting in supply chains, verifying its effectiveness in practical applications. \cite{Zhang2024} proposed the additive matrix autoregressive model. \cite{Tsay2024} provided a review of recent research on matrix time series analysis and proposed the matrix autoregressive moving average model, further improving the model system for matrix time series. 

In addition to the continuous multivariate time series data discussed above, non-negative integer-valued multivariate time series data are also very common in real-world applications. For example, the per-second frequency of sensitive words recorded from different IP addresses provides an important resource for monitoring the evolution of online public opinion. Numerous studies have shown that such word frequency statistics can effectively capture trends in online public sentiment. Similarly, daily reports of suspected cases, confirmed cases, and deaths across multiple countries constitute a crucial epidemic transmission matrix, which played an indispensable role in the prevention and control of the COVID-19 pandemic. Furthermore, monthly crime counts disaggregated by crime type and city have been shown to facilitate the analysis of spatiotemporal patterns in criminal activity and assist in criminal investigations. However, traditional continuous-valued multivariate time series models are often unsuitable for such data. These models not only fail to provide the integer-valued forecasts required for practical interpretation and decision-making, but also lack inherent interpretability when applied to discrete data. As noted in the literature, autoregressive moving average models, which assume real-valued outcomes and normally distributed errors, are generally inappropriate for count time series data. Therefore, there is an urgent need to establish a reasonable statistical framework for matrix-variate integer-valued time series data, a direction that has received increasing attention in recent methodological research.

In the modeling of integer-valued time series, \cite{Steutel1979} introduced the binomial thinning operator, laying a crucial foundation for this field. 
\cite{Alosh1987} subsequently proposed the first-order integer-valued autoregressive INAR(1) model based on this operator, marking the formal inception of integer-valued time series research. As the field has progressed, integer-valued time series models have been extended in various directions. The thinning operator, as a classic tool for count data modeling, possesses the core advantage of generating and preserving the integer-valued nature of the data. Compared to methods that treat count data as continuous, models based on the thinning operator construct a discrete stochastic process framework, endowing model parameters with intuitive probabilistic interpretations and thereby enhancing both interpretability and reliability (\cite{Wei2008}).
Building upon this foundation, \cite{Aknouche2024} further established the multiplicative integer-valued GARCH model, extending the volatility modeling capability of integer-valued time series. Collectively, these studies have enriched the methodological system of integer-valued time series and laid an important foundation for research on matrix-valued integer-valued time series.

Overall, modeling matrix integer-valued time series involves at least two key issues that need to be addressed. The first issue is the modeling mechanism for matrix-valued count data. Traditional thinning operators are only applicable to point-valued or vector data. Therefore, it is necessary first to define new thinning operations specifically designed for matrix-valued counts. The second issue is how to control the number of model parameters; too many parameters can lead to reduced estimation accuracy. These issues also point the way for subsequent research. Although research on integer-valued time series models has matured, their extension to the matrix case is still in its infancy. It was not until 2024 that \cite{Xu2024} systematically proposed a matrix integer-valued autoregressive model, establishing a complete modeling framework for the first time. This model defines a matrix form of the binomial thinning operator, replacing element-wise operations with binomial thinning operations while retaining the matrix multiplication structure, thereby balancing discrete and structural characteristics. 

However, as \cite{Ristic2009} pointed out, the Poisson distribution is not always suitable for modeling integer-valued time series, due to its inherent property that the mean equals the variance, a condition frequently violated in real-world data. In contrast, the negative binomial thinning operator, associated with geometric counting sequences, can better accommodate overdispersed data features and thus demonstrates significant advantages. Time series models based on the negative binomial thinning operator have important applications in fields such as reliability theory, medicine, and meteorology, with typical examples including the number of equipment awaiting maintenance, the number of congenital malformation cases, consecutive dry and wet days, and the occurrence of major global earthquakes. In this paper, we develop an autoregressive model for matrix-valued integer-valued data based on the negative binomial thinning operator, aiming to provide an effective methodological tool for the analysis of overdispersed matrix time series.

The remainder of this paper is organized as follows. In Section \ref{sec:2}, the definition and properties of the matrix negative binomial thinning operator are presented. In Section \ref{sec:3}, we introduce the matrix integer-valued autoregressive model based on the negative binomial thinning operator and discuss some of its probabilistic properties. In Section \ref{sec:4}, considers the projection estimation method and presents the asymptotic theory for the estimators. In Section \ref{sec:5}, the iterative conditional least squares estimation and its associated asymptotic theory are introduced. Section \ref{sec:6} discusses the order determination problem for the proposed model. Section \ref{sec:7} presents some numerical studies. Section \ref{sec:8} provides a practical application example.

\section{Matricial negative thinning operators}\label{sec:2}

Firstly, to address the limitations where the Poisson distribution is not always suitable for modeling integer-valued time series, we introduce the definitions of the left negative binomial thinning operator ``$\bm{A}{\ast}_L$'' and the right negative binomial thinning operator ``${\ast}_R\bm{B}$'' in this study.
Let $\bm{A}=(a_{i,j})_{m\times m}$ be a square matrix of order $m$  with elements $a_{i,j}\in [0, 1)$ and  $\bm{B}=(b_{i,j})_{n\times n}$ be a square matrix of order $n$  with elements $b_{i,j}\in [0, 1)$.
Denote by $\bm{Y}$  an  integer-valued matrix variate, i.e.,
$\bm{Y}=(y_{i,j})_{m\times n}$  with $y_{i,j} \in \mathbb{N}_0$.
Then, we define the $(i,j)$ elements of $\bm{A}{\ast}_L\bm{Y}$ and $\bm{Y}{\ast}_R\bm{B}$ by
\begin{equation}\label{ast}
	(\bm{A}{\ast}_L\bm{Y})_{i,j}=\sum_{k=1}^m a_{i,k}\ast y_{k,j}, ~(\bm{Y}{\ast}_R\bm{B})_{i,j}=\sum_{k=1}^n b_{k,i}\ast y_{j,k},
\end{equation}
where
\begin{equation}\label{star operator}
	a_{i,k}\ast y_{k,j}=
	\begin{cases}
		\sum_{s=1}^{y_{k,j}} W_{a_{i,k}}^{(s)},\mbox{~if~}y_{k,j}>0,\\
		0,\hspace{0.6in}\mbox{~if~}y_{k,j}=0,
	\end{cases}
	~b_{k,i}\ast y_{j,k}=
	\begin{cases}
		\sum_{s=1}^{y_{j,k}} W_{b_{k,i}}^{(s)},\mbox{~if~}y_{j,k}>0,\\
		0,\hspace{0.6in}\mbox{~if~}y_{j,k}=0,
	\end{cases}
\end{equation}
$\{W_{ai,k}^{(s)}\}$ and $\{W_{bk,i}^{(s)}\}$ are independent and identically distributed (i.i.d.) random variables following a geometric distribution with success probability$\dfrac{a_{i,k}}{1+a_{i,k}}$ and $\dfrac{b_{k,i}}{1+b_{k,i}}$.

Based on (\ref{ast}) and (\ref{star operator}), the left and right matrix thinning operators follow a rule analogous to that of standard matrix left and right multiplication, except that scalar multiplication is replaced by the binomial thinning operator in the element-wise computation.
Then, we give some useful properties of the matricial negative thinning operators ``$\bm{A}{\ast}_L$'' and ``${\ast}_R\bm{B}$'' in the following Proposition \ref{pro1}.
\begin{proposition}\label{pro1}
Consider $\bm{Y}=(y_{i,j})_{m\times n}$ and $\bm{Z}=(z_{i,j})_{m\times n}$ as two integer-valued matrix variates, where $y_{i,j},~z_{i,j} \in \mathbb{N}_0$.
For the matricial negative thinning operators defined in {\rm (\ref{ast})} and {\rm (\ref{star operator})}, the subsequent properties are satisfied.

\begin{enumerate}[(\text{\bf P}1)]
 \setlength{\itemsep}{0pt}
 \setlength{\parskip}{0pt}
 \setlength{\parsep}{0pt}
  \item Let $\bm{0}$ be a  matrix with all elements equal to 0, $\bm{I}$ be an identity matrix. Then,\\
  $(i)~\bm{0}_{m\times m}{\ast}_L\bm{Y}=\bm{0}_{m\times n},~\bm{Y}{\ast}_R\bm{0}_{n\times n}=\bm{0}_{m\times n}$; \\
  $(ii)~\bm{I}_{m\times 
  m}{\ast}_L\bm{Y}=\bm{Y},~\bm{Y}{\ast}_R\bm{I}_{n\times 
  n}=\bm{Y},~\mbox{and}$ \\
  $(iii)~\bm{0}_{m\times m}{\ast}_L\bm{Y}{\ast}_R\bm{I}_{n\times n}=\bm{0}_{m\times n},~\bm{I}_{m\times m}{\ast}_L\bm{Y}{\ast}_R\bm{0}_{n\times n}=\bm{0}_{m\times n}$.
  \item
  For a constant $c$ with $0 \leq c \leq 1$, we have\\
  $
 (iv)~ \bm{A}{ \ast}_{L} \bm{Y}{ \ast}_{R}~c=(c\bm{A}){ \ast}_{L} \bm{Y}, c~{ \ast}_{L} \bm{Y}{ \ast}_{R}\bm{B}^{\top}=\bm{Y}{ \ast}_{R}(c\bm{B})^{\top}, and
 ~(v)~(c \bm{A}){ \ast}_{L} \bm{Y}{ \ast}_{R}\bm{B}^{\top}=\bm{A}{ \ast}_{L} \bm{Y}{ \ast}_{R}( c \bm{B})^{\top}.
  $
  \item  The operators ``$\bm{A}{\ast}_L$'' and ``${\ast}_R\bm{B}$'' adhere to the following operational rules: \\
     $(vi)~(\bm{A}{\ast}_L\bm{Y}){\ast}_R\bm{B}=\bm{A}{\ast}_L(\bm{Y}{\ast}_R\bm{B})$;\\
     $(vii)~(\bm{A}{\ast}_L\bm{Y})^{\top}=\bm{Y}^{\top}{\ast}_R\bm{A}^{\top}$, $(\bm{Y}{\ast}_R\bm{B})^{\top}=\bm{B}^{\top}{\ast}_L\bm{Y}^{\top}$, and $(\bm{A}{\ast}_L\bm{Y}{\ast}_R\bm{B})^{\top}=\bm{B}^{\top}{\ast}_L\bm{Y}^{\top}{\ast}_R\bm{A}^{\top}$;\\
     $(viii)~\bm{A}{\ast}_L(\bm{Y}+\bm{Z})=\bm{A}{\ast}_L\bm{Y}+\bm{A}{\ast}_L\bm{Z}$,  and
     $(\bm{Y}+\bm{Z}){\ast}_R\bm{B}=\bm{Y}{\ast}_R\bm{B}+\bm{Z}{\ast}_R\bm{B}$.
  \item The linearity of mathematical expectations holds, namely\\
  $(ix)~E(\bm{A}{\ast}_L\bm{Y})=\bm{A}E(\bm{Y}), ~E(\bm{Y}{\ast}_R\bm{B})=E(\bm{Y})\bm{B}$, and $E(\bm{A}{\ast}_L\bm{Y}{\ast}_R\bm{B})=\bm{A}E(\bm{Y})\bm{B}$.

\end{enumerate}
\end{proposition}

The preceding properties are useful for gaining insight into the matricial negative thinning operators and for deriving results related to the proposed matrix‑variate integer-valued autoregressive model. In what follows, we denote  $\bm{I}_{m\times m}$ simply by $\bm{I}_m$ when no confusion arises. Likewise, for simplicity, we write $\bm{0}_m$ to represent an $m$-dimensional square matrix in which every entry is zero, instead of using the notation $\bm{0}_{m\times m}$.
\section{Matrix integer-valued autoregressive process with negative binomial thinning (MAT-NB-INAR)}\label{sec:3}
\subsection{Definition of MAT-NB-INAR($p$) process}

We introduce the MAT-NB-INAR($p$) process in the following recursive equation:
\begin{equation}\label{Mymodelp}
	\bm{Y}_t = \bm{A}_1{ \ast}_{L} \bm{Y}_{t-1}{ \ast}_{R}\bm{B}_1^{\top}+\cdots+\bm{A}_p{ \ast}_{L} \bm{Y}_{t-p}{ \ast}_{R}\bm{B}_p^{\top}+ \bm{\mathcal{E}}_{t} ,~t \in \mathbb{Z},
\end{equation}
where
\begin{enumerate}[(i)]
  \setlength{\itemsep}{0pt}
  \setlength{\parskip}{0pt}
  \setlength{\parsep}{0pt}
\item $\bm{Y}_t$ denotes the $m\times n$  matrix-variate integer-valued observation at time $t,~t\in\{1,...,T\}$;
\item $\bm{A}_l=(a_{i,j}^{(l)})_{m \times m}$ is $m\times m$ autoregressive coefficient matrices whose entries satisfyis $a_{i,j}^{(l)} \in [0,1)$, and
  $\bm{B}_l=(b_{s,k}^{(l)})_{n\times n}$ is the $n\times n$ autoregressive coefficient matrices with entries $(b_{s,k}^{(l)})_{n\times n}\in [0,1)$, where $i,j\in\{1,...,m\}$, $s,k\in\{1,...,n\}$, $l\in\{1,...,p\}$; to ensure model identifiability, we normalize $\bm{A}_l$ such that ${\Vert\bm{A}_l\Vert}_F=1$ for {\rm($l\in \{1,...,p\}$)}.
\item The symbols ``$\bm{A}_l{\ast}_L$'' and ''${\ast}_R\bm{B}_l$'' ($l\in\{1,...,p\}$) represent the matricial negative thinning operators
	  introduced in (\ref{ast}) and (\ref{star operator});
\item $\{\bm{\mathcal{E}}_t\}$  denotes a sequence of i.i.d. $\mathbb{N}^{m\times n}$-valued random error matrix
      following some discrete matrix-valued distribution with mean matirx $\bm{\Lambda}$, and satisfying
	  $E(\bm{\mathcal{E}}_t\bm{\mathcal{E}}_t^{\top}) < \infty$; 
\item For any fixed $t$, $\bm{\mathcal{E}}_t$ is taken to be independent of the counting series involved in $\bm{A}_l{\ast}_L\bm{Y}_{t-l}{\ast}_R\bm{B}_l^{\top}$ ($l\in\{1,...,p\}$) and $\bm{Y}_{t-s}$ for all $s\geq 1$; moreover, $\bm{A}_l{\ast}_L\bm{Y}_{t-l}{\ast}_R\bm{B}_l^{\top}$ ($l\in\{1,...,p\}$) and  $\bm{Y}_{t-s}$ ($s\geq 1$) are independent of each other.
\end{enumerate}

We define ${\rm vec}(\cdot)$ as the vectorization operator that stacks the columns of a matrix, and subsequently rewrite the  MAT-NB-INAR($p$) model in (\ref{Mymodelp}) as:
\begin{equation}\label{Mymodel1}
{\rm vec}(\bm{Y}_t) =(\bm{B}_1\otimes \bm{A}_1){ \ast}_{L} {\rm vec}(\bm{Y}_{t-1})+\cdots+(\bm{B}_p\otimes \bm{A}_p){ \ast}_{L} {\rm vec}(\bm{Y}_{t-p}) +{\rm vec}(\bm{\mathcal{E}}_{t}) ,~t \in \mathbb{Z},
\end{equation}
where ``$\otimes$'' is Kronecker product (\cite{Liu2008}).
If we further define this in terms of $\bm{\Phi}_l=\bm{B}_l \otimes \bm{A}_l$ $(l\in\{1,...,p\})$, then the model can be rewritten as:
\begin{align}\label{vector_model}
{\rm vec}(\bm{Y}_t) =\bm\Phi_1{\ast}_{L} {\rm vec}(\bm{Y}_{t-1})+\cdots+ \bm\Phi_p{\ast}_{L} {\rm vec}(\bm{Y}_{t-p})+{\rm vec}(\bm{\mathcal{E}}_{t}) ,~t \in \mathbb{Z}.
\end{align}

Equation (\ref{vector_model}) essentially corresponds to the multivariate generalized INAR (MGINAR($p$)) model proposed by \cite{Latour1997}, which involves $p(mn)^2+nm$ parameters to be estimated, whereas model (\ref{Mymodelp}) contains only $p(m^2+n^2)+mn$ parameters. We observe that $p(mn)^2+nm$  is significantly larger than $p(m^2+n^2)+mn$. Moreover, model (\ref{Mymodelp}) offers stronger interpretability by preserving the matrix structure of the observed data. The coefficient matrices  $\bm{A}_1,...,\bm{A}_p$ and $\bm{B}_1,...,\bm{B}_p$  in model (\ref{Mymodelp}) play a crucial role. Specifically, the left matrices $\bm{A}_1,...,\bm{A}_p$ and the right matrices $\bm{B}_1,...,\bm{B}_p$ capture row-wise and column-wise dependencies, respectively. Thus, model (\ref{Mymodelp}) not only accounts for the interactive correlations between the rows and columns of the matrix-variate observations but also characterizes the autocorrelation structure of the matrix time series over time. If model (\ref{vector_model}) are used for modeling, the inherent structure of the matrix data would be disrupted, leading to substantial information loss.

To illustrate the necessity of applying negative binomial thinning operators for modeling, we first analyze a restricted scenario involving only a single negative binomial thinning operator, beginning with the case where the left thinning operator is retained in isolation. If we take $\bm{B}_1=\cdots=\bm{B}_p=\bm{I}_n$, then model (\ref{Mymodelp}) simplifies to:
\begin{equation}\label{rowmodel}
	\bm{Y}_t=\bm{A}_1{\ast}_L\bm{Y}_{t-1}+\cdots+\bm{A}_p{\ast}_L\bm{Y}_{t-p}+\bm{\mathcal{E}}_{t}.
\end{equation}
We proceed to interpret and analyze the model from the row-wise perspective. Let $\bm{Y}_{i,\cdot,t}:=(y_{i1,t},y_{i2,t},...,y_{in,t})$ denote the $i$th row vector of the matrix $\bm{Y}_t$, and similarly, let the $i$th row vector of the corresponding disturbance matrix $\bm{\mathcal{E}}_t$ be denoted as $\bm{\mathcal{E}}_{i,\cdot,t}:=(e_{i1,t},e_{i2,t},...,e_{in,t})$. 
By decomposing equation (\ref{rowmodel}) along the row dimension and applying a transpose transformation, we obtain the autoregressive expression in the row dimension:
$$
	\bm{Y}_{i,\cdot,t}^{\top}= \sum_{l=1}^p \left(a_{i,1}^{(l)} \ast \bm{Y}_{1,\cdot,t-l}^{\top} + a_{i,2}^{(l)} \ast \bm{Y}_{2,\cdot,t-l}^{\top} +
	\cdots + a_{i,m}^{(l)} \ast \bm{Y}_{m,\cdot,t-l}^{\top}\right) + \bm{\mathcal{E}}_{i,\cdot,t}^{\top},~i\in\{1,...,m\},
$$
where $a_{i,s}^{(l)}\ast\bm{Y}_{i,\cdot,t-l}=(a_{i,s}^{(l)}\ast y_{i1,t-l},...,a_{i,s}^{(l)}\ast y_{in,t-l})$, $s\in\{1,...­,m\}$.
This decomposition shows that model (8) splits into $m$ sub-models, each taking an autoregressive form built from the rows of  $\bm{Y}_{t-1},\bm{Y}_{t-2},...,\bm{Y}_{t-p}$. Hence, when only the left thinning operation is adopted, the model accounts solely for row-wise dependencies, neglecting the column-wise interactions.

Conversely, if we let $\bm{A}_1=\cdots=\bm{A}_p=\bm{I}_m$, model (\ref{Mymodelp}) reduces to:
\begin{equation}\label{columnmodel}
	\bm{Y}_t=\bm{Y}_{t-1}{\ast}_R\bm{B}_1^{\top}+\cdots+\bm{Y}_{t-p}{\ast}_R\bm{B}_p^{\top}+\bm{\mathcal{E}}_{t}.
\end{equation}
Let $\bm{Y}_{\cdot,j,t}:=(y_{1j,t},y_{2j,t},...,y_{mj,t})^{\top}$ be the $j$th column of $\bm{Y}_t$ and similarly define  $\bm{\mathcal{E}}_{\cdot,j,t}:=(e_{1j,t},e_{2j,t},...,e_{mj,t})^{\top}$ as the $j$th column of $\bm{\mathcal{E}}_t$.
Then we obtain
\begin{equation}\label{columnsubmodel}
	\bm{Y}_{\cdot,j,t}=\sum_{l=1}^p \left(b_{j,1}^{(l)} \ast \bm{Y}_{\cdot,1,t-l} + b_{j,2}^{(l)} \ast \bm{Y}_{\cdot,2,t-l} +\cdots+ b_{j,n}^{(l)} \ast \bm{Y}_{\cdot,n,t-l}\right)+\bm{\mathcal{E}}_{\cdot,j,t},~j=1,...,n,
\end{equation}
when $b_{j,k}^{(l)}\ast\bm{Y}_{\cdot,j,t-l}=(b_{j,k}^{(l)}\ast y_{1j,t-l},...,b_{j,k}^{(l)}\ast y_{mj,t-l})^{\top}$, for $k\in\{1,...,n\}$. From (\ref{columnsubmodel}) it is evident that model (\ref{columnmodel}) decouples into $n$ sub-models, each of which follows an autoregressive structure constructed from the columns of the matrices $\bm{Y}_{t-1},\bm{Y}_{t-2},...,\bm{Y}_{t-p}$. Consequently, when only the right thinning operation is employed, the resulting model focuses exclusively on modeling the columns of the matrix-variate observations, while the information embedded in the rows is disregarded.

{\subsection{Probabilistic properties}\label{model_properties}

In the following proposition, we state the existence, strict stationarity and ergodicity of the MAT-NB-INAR($p$) process defined in (\ref{Mymodelp}). We define
$$
{\bm{\mathcal{A}}}=\left(\begin{array}{lllll}
\bm{B}_1\otimes \bm{A}_1&\bm{B}_2\otimes \bm{A}_2&\cdots&\bm{B}_{p-1}\otimes \bm{A}_{p-1}&\bm{B}_p\otimes \bm{A}_p\\
\mathbf{I}_{mn}&\mathbf{0}_{mn}&\cdots&\mathbf{0}_{mn}&\mathbf{0}_{mn}\\
\vdots&\vdots&\ddots&\vdots&\vdots\\
\mathbf{0}_{mn}&\mathbf{0}_{mn}&\cdots&\mathbf{I}_{mn}&\mathbf{0}_{mn}
\end{array}\right).
$$
Thus, we have the following proposition. 
\begin{proposition}\label{stationary}
If  $\rho({\bm{\mathcal{A}}})<1$, then the MAT-NB-INAR($p$) process defined in $(\ref{Mymodelp})$ is stationary and causal, where for any square matrix,  $\rho(\cdot)$  denotes its spectral radius, i.e., the maximum modulus of the (complex) eigenvalues of this matrix.
\end{proposition}

Under the conditions of Proposition \ref{stationary}, derivations for the conditional expectation, unconditional expectation and variance-covariance matrix closely follow the framework developed by \cite{Xu2024}. This proposition is fundamental to the theoretical development of the MAT-NB-INAR($p$) model, as it guarantees the stability of the process and forms the basis for deriving the asymptotic properties of the parameter estimators discussed in Section \ref{sec:4}.
Unlike the standard MAT-INAR(p) process, the use of negative binomial thinning does not alter the stationarity condition, as the thinning operation preserves the linear autoregressive structure in the conditional mean. Consequently, the stability condition depends solely on the autoregressive coefficient matrices $\bm{A}_l$ and $\bm{B}_l$.

When the condition of Proposition \ref{stationary} is satisfied, we study the properties of the conditional moments and moments of the MAT-NB-INAR($p$) process. We first derive the conditional expectation and expectation for the MAT-NB-INAR($p$) process as
\begin{align}
    E(\bm{Y}_t|\bm{Y}_{t-1}, \ldots, \bm{Y}_{t-p}) 
    &= \sum_{l=1}^p \bm{A}_l \bm{Y}_{t-l} \bm{B}_l^{\top} + \bm{\Lambda}, \quad t \in \mathbb{Z}. \label{e1} \\
    \intertext{Taking unconditional expectations on both sides of \eqref{e1} and solving for $\operatorname{vec}(\bm{\mu}_Y) := E(\operatorname{vec}(\bm{Y}_t))$, we obtain}
    \operatorname{vec}(\bm{\mu}_Y) 
    &= \left( \bm{I}_{mn} - \sum_{l=1}^p \bm{B}_l \otimes \bm{A}_l \right)^{-1} \operatorname{vec}(\bm{\Lambda}), \quad t \in \mathbb{Z}. \label{e2}
\end{align}
Equation \eqref{e1} can be obtained from property ({\bf P}4) in Section \ref{sec:2}. Based on \cite{Latour1997}, equation \eqref{e2} can be readily proven for model \eqref{Mymodel1}.

We now derive the variance-covariance matrix of $\bm{Y}_t$, which is of particular importance for matrix-valued time series models as it reveals the cross-correlation structure along both the row and column dimensions, thereby enabling a simultaneous examination of the interdependencies among all component series. Let $\bm{\Gamma}_h := {\rm Cov}({\rm vec}(\bm{Y}_t),{\rm vec}(\bm{Y}_t-h))$ denote the autocovariance matrix at lag $h$.
To characterize the correlation structure of $\bm{Y}_t$ in the row and column directions, we introduce the transformation matrix 
$\mathcal{T} = (\tau_{i,j})_{mn \times mn}$, where $\tau_{i,j} = 1$ if and only if  $i \in \{sm + 1, sm + 2, \ldots, (s+1)m\}$, $j \in \{(s+1) + (i-sm-1)n\}$ with $s \in \{0, 1, \ldots, n-1\}$, and $\tau_{i,j} = 0$ otherwise. This transformation matrix rearranges the Kronecker product-based covariance matrix into a block structure, facilitating the identification of correlations between rows and between columns.

Specifically, the Kronecker-based covariance matrix is defined as
$$
\bm{\Gamma}_0^{\otimes}=E\left((\bm{Y}_t-\bm{\mu}_Y)\otimes(\bm{Y}_t-\bm{\mu}_Y)^{\top}\right),~t \in \mathbb{Z}.
$$
Upon applying the transformation matrix $\mathcal{T}$, the column-wise and row-wise covariance matrices are respectively obtained as
$$
\bm{\Sigma}_0^{c}:=\bm{\mathcal{T}}\bm{\Gamma}_0^{\otimes},\bm{\Sigma}_0^{r}:=\bm{\Gamma}_0^{\otimes}\bm{\mathcal{T}}.
$$
$\bm{\Sigma}_0^{c}$ is rearranged into an $n\times n$ block matrix consisting of $m \times m$ submatrices, where the diagonal blocks represent the variances of individual columns and the off-diagonal blocks represent the covariances between different columns. Similarly, $\bm{\Sigma}_0^{r}$ is rearranged into an $m \times m$ block matrix consisting of $n\times n$ submatrices, with diagonal blocks corresponding to the variances of individual rows and off-diagonal blocks corresponding to the covariances between different rows. This decomposition enables us to separately examine the contemporaneous correlations between column direction and between row direction. Extending the above definitions to the lag-$h$ case, let
$$
\bm{\Gamma}_h^{\otimes}=E\left((\bm{Y}_{t+h}-\bm{\mu}_Y)\otimes(\bm{Y}_t-\bm{\mu}_Y)^{\top}\right),~t \in \mathbb{Z}.
$$
Then the lag-$h$ column-wise and row-wise autocovariance matrices are given by
$$
\bm{\Sigma}_h^{c}:=\bm{\mathcal{T}}\bm{\Gamma}_h^{\otimes},\bm{\Sigma}_h^{r}:=\bm{\Gamma}_h^{\otimes}\bm{\mathcal{T}}.
$$
It is worth noting that for a matrix-variate integer-valued time series $\{\bm{Y}_t\}$, $\bm{\Gamma}_h^{\otimes}\neq\bm{\Gamma}_{-h}^{\otimes}$ (\cite{Samadi2014}).
The relationship between $\bm{\Gamma}_h^{\otimes}$ and $\bm{\Gamma}_{-h}^{\otimes}$ is given by
${\bm{\Gamma}_h^{\otimes}}=(\bm{\mathcal{T}}\bm{\Gamma}_{-h}^{\otimes}\bm{\mathcal{T}})^{\top}$.

\section{Projection estimation}\label{sec:4}

\indent Suppose $\{\bm{Y}_t\}_{t=1}^T$ is a sequence of matrix-valued observations generated from the MAT-NB-INAR($p$) process. Define the parameter space for ($\bm{A}_1,\bm{B}_1,...,\bm{A}_p, \bm{B}_p$, $\bm{\Lambda}$) as
$$
\bm{\Theta}=
\{
\underbrace{{[0,1]}^{m\times m}\times{[0,1]}^{n\times n}\times \cdots \times{[0,1]}^{m\times m}\times{[0,1]}^{n\times n}}_{p~\text{times the product of}~{[0,1]}^{m\times m}\times{[0,1]}^{n\times n}}\times(0,\infty)^{m\times n}
\},
$$
where the first $2p$ components correspond to the autoregressive coefficient matrices and the final component corresponds to the mean matrix of the error term. In what follows, we investigate the projection method for estimating the $\bm{A}_1,\bm{B}_1,...,\bm{A}_p, \bm{B}_p$ and $\bm{\Lambda}$, along with the asymptotic properties of the resulting estimators.

\subsection{Projection method}

In this section, we apply a two-stage projection approach to estimate the parameters $\bm{A}_1,\bm{B}_1,...,\bm{A}_p, \bm{B}_p$ and $\bm{\Lambda}$. The first stage focuses on the vectorized representation of the model and yields preliminary estimates via conditional least squares. The second stage recovers the original matrix coefficient matrices by solving a nearest Kronecker product (NKP) problem.

The first stage is vectorized estimation. Recall that $\bm{\Phi}_l=\bm{B}_l \otimes \bm{A}_l~(l\in\{1,...,p\})$, we obtain the conditional least squares (CLS) estimates $\widehat{\bm{ \Phi}}_1,...,\widehat{\bm{ \Phi}}_p$ and $\widehat{\bm\Lambda}$ for $\bm{\Phi}_1,...,\bm{\Phi}_p$ and $\bm{\Lambda}$.
Specifically, define $\bm{\Psi}^{\top}:=(\bm\Phi_1,...,\bm\Phi_p, {\rm vec}(\bm\Lambda))$,
$\bm{X}_{t}^{\top}:=({\rm vec}(\bm{Y}_{t})^{\top},...,{\rm vec}(\bm{Y}_{t-p+1})^{\top}, {1})$, and let
${\bm{\mathcal X}}=(\bm{X}_{p},...,\bm{X}_{T-1})^{\top}$, and $\bm{\mathcal Y}=({\rm vec}(\bm{Y}_{p+1}),...,{\rm vec}(\bm{Y}_{T}))^{\top}$. The CLS function for $\bm\Psi$ is
$$
Q(\bm\Psi):=\left(\bm{\mathcal Y}-G(\bm\Psi)\right)^{\top}\left(\bm{\mathcal Y}-G(\bm\Psi)\right),\notag
$$
with $G(\bm\Psi):=E(\bm{\mathcal Y}|{\rm vec}(\bm{Y}_{t-1}),...,{\rm vec}(\bm{Y}_{t-p}))={\bm{\mathcal X}}\bm\Psi$. Minimizing $Q(\bm\Psi)$ over the parameter space $\bm{\Theta}^{\ast} \subseteq \mathbb{R}^{(pmn+1)\times (mn)}$ yields the CLS estimator
\begin{equation}\label{clspsi}
\widehat{\bm{\Psi}}=(\bm{\mathcal  X}^{\top}\bm{\mathcal X})^{-1}\bm{\mathcal  X}^{\top}\bm{\mathcal Y}.
\end{equation}
From $\widehat{\bm{\Psi}}$ we directly obtain  $\widehat{\bm \Phi}_1,...,\widehat{\bm \Phi}_p$ and $\widehat{\bm\Lambda}$, and set $\widehat{\bm\Lambda}_{PROJ}=\widehat{\bm\Lambda}$.

The second stage is the nearest Kronecker product decomposition. For each $l\in \{1,...,p\}$, we aim to find $\bm{A}_l$ and $\bm{B}_l$ that best approximate $\widehat{\bm \Phi}_l$ in the Frobenius norm subject to a NKP structure:
\begin{equation}\label{NKP}
(\widehat{\bm{A}}_l,\widehat{\bm{B}}_l)=\mathop{\arg\min}\limits_{\bm{A}_l, \bm{B}_l}{\left\Vert \widehat{\bm\Phi}_l-\bm{B}_l\otimes \bm{A}_l\right\Vert}_F^2,~l\in \{1,...,p\},
\end{equation}
This is a standard nearest NKP problem (\cite{Loan2000}).

To solve (\ref{NKP}), we employ the re-arrangement operator $g :\mathbb{R}^{mn}\times \mathbb{R}^{mn}\rightarrow \mathbb{R}^{m^2}\times \mathbb{R}^{n^2},$ introduced by \cite{Chen2021}. The operator rearranges its entries so that $g(\bm{B}_l\otimes\bm{A}_l)={\rm vec}(\bm{A}_l){\rm vec}(\bm{B}_l)^{\top}$. 
Operationally, we first partition $\bm{\Phi}_l$ according to the NKP structure:
$$
\bm{\Phi}_l=\left(\begin{array}{ccc}
\bm{\Phi}_{1,1}^{(l)} &\cdots & \bm{\Phi}_{1,n}^{(l)}\\
\vdots  & \ddots  &\vdots \\
\bm{\Phi}_{n,1}^{(l)} &... & \bm{\Phi}_{n,n}^{(l)}\\
\end{array}\right),
\notag
$$
where each block $\bm{\Phi}_{i,j}^{(l)}=b_{i,j}^{(l)}\bm{A}_{l}$ is of size $m\times m$.
The rearrangement then gives
$$
\widetilde{\bm\Phi}:=g(\widehat{\bm\Phi})=\left(
{\rm vec}(\widehat{\bm\Phi}_{1,1}^{(l)}),....,{\rm vec}(\widehat{\bm \Phi}_{n,1}^{(l)}),...,{\rm vec}(\widehat{\bm \Phi}_{1,n}^{(l)}),...,{\rm vec}(\widehat{\bm \Phi}_{n,n}^{(l)})\right).
$$
For two $m\times n$ matrices $\bm{C}$ and $\bm{D}$, the operator ${g}$ is linear and preserves the Frobenius norm, i.e.:

(i) $g(\bm{C}+\bm{D})=g(\bm{C})+g(\bm{D})$, and 

(ii) ${\Vert g(\bm{C})\Vert}_F={\Vert\bm{C}\Vert}_F$.\\
Consequently, problem (\ref{NKP}) is equivalent to
\begin{align}
(\widehat{\bm{A}}_l,\widehat{\bm{B}}_l)=\mathop{\arg\min}\limits_{\bm{A}_l,\bm{B}_l}{\Vert \widehat{\bm\Phi}_l-\bm{B}_l\otimes\bm{A}_l\Vert}_F^2&=\mathop{\arg\min}\limits_{\bm{A}_l,\bm{B}_l}{\Vert g(\widehat{\bm\Phi}_l)-g(\bm{B}_l\otimes\bm{A}_l)\Vert}_F^2\notag\\
&=\mathop{\arg\min}\limits_{\bm{A}_l,\bm{B}_l}{\Vert \widetilde{\bm\Phi}_l-{\rm vec}(\bm{A}_l){\rm vec}(\bm{B}_l)^{\top}\Vert}_F^2,~l\in \{1,...,p\}.
\notag
\end{align}
Following the standard NKP solution (\cite{Loan2000}), we compute the singular value decomposition (SVD) of $\widetilde{\bm\Phi}_l$. Let $d_{1}^{(l)}$ be the largest singular value of $\widetilde{\bm\Phi}_l$, and let $\bm{u}_{1}^{(l)}$ and $\bm{v}_{1}^{(l)}$ be the corresponding left and right singular vectors. Then the solution is given by
$$
{\rm vec}({\bm{{A}}}_{l})=\sqrt{d_{1}^{(l)}}\bm{u}_{1}^{(l)},~~{\rm vec}({\bm{{B}}}_{l})^{\top}=\sqrt{d_{1}^{(l)}}{\bm{v}_{1}^{(l)}}^{\top},~l\in\{1,...,p\}.
\notag
$$
The final estimates ${\widehat{\bm {A}}}_{l,PROJ}$ and ${\widehat{\bm {B}}}_{l,PROJ}$ are obtained by reshaping $\bm{u}_{1}^{(l)}$ and $d_1^{(l)}\bm{v}_1^{(l)}$ into matrices of appropriate dimensions. Notice that $\Vert\bm{{A}}_{l}\Vert_F=({d_{1}^{(l)}})^{1/2}$, which implies $\Vert\widehat{\bm{A}}_{l,PROJ}\Vert_F=1$ for $l\in\{1,...,p\}$, thereby ensuring the identifiability of model (\ref{Mymodelp}).

\subsection{The asymptotic property of projection estimator}
We now establish the asymptotic behavior of the projection (PROJ) estimators introduced in the previous section. The analysis proceeds in two steps: first, we show that the vectorized coefficient estimators $\bm{\hat \Phi}_l~(l \in \{1,...,p\})$ are asymptotically normal; then, we derive the asymptotic distributions of the reconstructed matrix estimators ${\widehat{\bm{A}}}_{l,PROJ}$,${\widehat{\bm {B}}}_{l,PROJ}$ and $\widehat{\bm \Lambda}_{PROJ}$.

\begin{theorem}
\label{thm1}
Let $\{\bm{Y}_t\}$ be a MAT-NB-INAR($p$) process satisfying {$\rho({\bm{\mathcal{A}}})<1$}, and assume that for each $l\in\{1,...,p\}$, the matrices $\bm{A}_l,~\bm{B}_l,~\bm{\Lambda}$, and the covariance matrix of $\bm{\mathcal{E}}_{t}$ are nonsingular. Then the estimator $\bm{\hat \Phi}_l$ obtained from (\ref{clspsi}) is strongly consistent and asymptotically normal:
$$
\sqrt{T-p}~{\rm vec}(\widehat{\bm{\Phi}}_l^{\top}-\bm{\Phi}_l^{\top})\overset{L}{\longrightarrow}N(\bm{0},\bm{W}_l{\widetilde{\bm{\Sigma}}}\bm{W}_l^{\top}),~l\in\{1,...,p\}.
$$
We first define the selection matrix $\bm{W}_l$ as a block matrix with the following structure:
$\bm{W}_l=(\bm{0}_{m^2n^2},...,\bm{I}_{m^2n^2},...,\bm{0}_{m^2n^2},\bm{0}_{mn})$,
where the $l$th block entry is an identity matrix of order $(mn)^2$, and all other block entries are zero matrices with matching dimensions.
We then specify the asymptotic covariance matrix $\widetilde{\bm\Sigma}$ the NKP: $\widetilde{\bm\Sigma}= \bm{\Sigma}_{\mathcal  U}\otimes\bm{H}^{-1} $. It consists of two components: the first is related to the model disturbance term, let $\bm{{\mathcal{U}}}_t={\rm vec}(\bm{Y}_t)-{\rm vec}(\bm{{\Lambda}})-\sum_{l=1}^p{\bm {\Phi}}_l{\rm vec}(\bm{Y}_{t-l})$. The matrix $\bm{\Sigma}_{\mathcal{U}}:=E(\bm{\mathcal U}_{t}\bm{\mathcal U}_{t}^{\top})$ denotes the unconditional variance-covariance matrix of the disturbance $\bm{{\mathcal{U}}}_t$.
The second is matrix $\bm{H}$, which is defined as the second-moment matrix of the regressor $\bm{ X}_{t}$, that is $\bm{H}:=E(\bm{ X}_{t}\bm{ X}_{t}^{\top})$. We further establish $\widehat{\bm\Sigma}_{\mathcal{U}}=\sum_{t=p+1}^T\widehat{\bm {\mathcal{U}}}_t\widehat{\bm {\mathcal{U}}}_t^{\top}/(T-mn-p)$ converges a.s. to $\bm{\Sigma}_{\mathcal{U}}$, and $\widehat{\bm{ H}}_T=\sum_{t=p+1}^T\bm{X}_{t}\bm{X}_{t}^{\top}/(T-p)$ converges a.s. to $\bm{H}$.
\end{theorem}

The result follows directly from the standard theory of MGINAR($p$) models (\cite{Latour1997}).

\begin{theorem}
\label{thm2}
Under the conditions of \ref{thm1}, the PROJ estimators ${\widehat{\bm{A}}}_{l,PROJ}$ and ${\widehat{\bm {B}}}_{l,PROJ}$ are strongly consistent and asymptotically normal. Specifically,
$$
\sqrt{T-p}\left(\begin{array}{c}
{\rm vec}(\widehat{\bm{A}}_{l,PROJ}-\bm{A}_l) \\
{\rm vec}(\widehat{\bm{B}}_{l,PROJ}-\bm{B}_l)
\end{array}\right)\overset{L}{\longrightarrow}N\left(\bm{0},\bm{V}_0^{(l)}\bm{\Xi}_1^{(l)}{\bm{V}_0^{(l)}}^{\top}\right),~l\in\{1,...,p\},
$$
with
\begin{align}
&\bm{V}_{0}^{(l)}:=\left(\begin{array}{c}
{\Vert\bm{B}_l\Vert}_F^{-1} [{\bm{\beta}_l^{(1)}}^{\top}\otimes(\bm{I}_{m^2}-\bm{\alpha}_l\bm{\alpha}_l^{\top})]\nonumber\\
\bm{I}_{n^2}\otimes\bm{\alpha}_l^{\top}
\end{array}\right),\nonumber
\end{align}
where  $\bm{\alpha}_l:={\rm vec}(\bm{A}_l)$, $\bm{\beta}_l:={\rm vec}(\bm{B}_l)$, $\bm{\beta}_l^{(1)}:=\bm{\beta}_l/ {\Vert\bm{\beta}_l\Vert}$, and {$\bm{\Xi}_1^{(l)}=\widetilde{\bm{\mathcal{T}}}\bm{W}_l\widetilde{\bm\Sigma}^{\top}\bm{W}_l^{\top}\widetilde{\bm{\mathcal{T}}}$. $\widetilde{\bm{\mathcal{T}}}$ is an order of $mn^2$ block matrix, each block $\bm{T}_{i,j}$ forms an $m\times m$ matrix that $\bm{I}_m$ when the index conditions $i=ks+1+(j-sm-1)n$ for $kn-n\leq s\leq kn-1$ and $k=1,\ldots,n$, $j=sm+1,sm+2,...,(s+1)m$ and $s=0,1,...,(n^2-1)$ are all satisfied, and equals $\bm{0}_m$ otherwise.

Moreover, the asymptotic distribution of the NKP of the two estimated vectors is given by
$$
\sqrt{T-p}\left({\rm vec}(\widehat{\bm{B}}_{l,PROJ})\otimes{\rm vec}(\widehat{\bm{A}}_{l,PROJ})-{\rm vec}(\bm{B}_l)\otimes{\rm vec}(\bm{A}_l)\right)\overset{L}{\longrightarrow}N\left(\bm{0},\bm{V}_1^{(l)}\bm{\Xi}_1^{(l)}{\bm{V}_1^{(l)}}^{\top}\right),
$$
where
$\bm{V}_1^{(l)}:=(\bm{\beta}_l^{(1)}{\bm{\beta}_l^{(1)}}^{\top})\otimes\bm{I}_{m^2}+\bm{I}_{n^2}\otimes(\bm{\alpha}_l\bm{\alpha}_l^{\top})-(\bm{\beta}_l^{(1)}{\bm{\beta}_l^{(1)}}^{\top})\otimes(\bm{\alpha}_l\bm{\alpha}_l^{\top})$.
}
\end{theorem}

\begin{theorem}
Under the same conditions, the estimator $\widehat{\bm \Lambda}_{PROJ}$ is strongly consistent and asymptotically normal:
$$
\sqrt{T-p}~{\rm vec}(\widehat{\bm \Lambda}_{PROJ}^{\top}-\bm{ \Lambda}^{\top})\overset{L}{\longrightarrow}N\left(\bm{0},\bm{W}_{p+1}\widetilde{\bm\Sigma}\bm{W}_{p+1}^{\top}\right),\nonumber
$$
where $\bm{W}_{p+1}=(\bm{0}_{m^2n^2},...,\bm{0}_{m^2n^2},\bm{I}_{mn})$ selects the last block corresponding to the mean parameter.
\end{theorem}

Under the conditions of (\ref{thm1}), we can prove (\ref{thm2}) by generalizing the proof arguments of \cite{Chen2021}.
\section{Iterated conditional least squares estimation}\label{sec:5}
The projection method introduced in \ref{sec:4} requires an initial estimation of the high-dimensional matrix $\bm{\Psi}$. When the dimension  $(pmn+1)\times(mn) $ of $\bm{\Psi}$ is large, this becomes a high-dimensional or even ultra-high-dimensional estimation problem, and the resulting projection estimators may suffer from reduced accuracy. To address this issue, we develop an iterated conditional least squares (ICLS) procedure for estimating the parameters $\bm{A}_1,\bm{B}_1,...,\bm{A}_p,\bm{B}_p$ and $\bm{\Lambda}$. The ICLS algorithm takes the PROJ estimators as initial values and then iteratively refines them to achieve improved accuracy. The ICLS estimators $\widehat{\bm{A}}_{l,ICLS}$, $\widehat{\bm{B}}_{l,ICLS}~(l\in\{1,...,p\})$ and $\widehat{\bm{\Lambda}}_{ICLS}$ are defined as the solution to the following minimization problem:
$$
(\widehat{\bm{A}}_1, \widehat{\bm{B}}_1,...,\widehat{\bm{A}}_p,\widehat{\bm{B}}_p,\widehat{\bm{\Lambda}})=
\mathop{\arg\min}\limits_{\bm{A}_1,\bm{B}_1,...,\bm{A}_p,\bm{B}_p,\bm{\Lambda}}\sum_{t=p+1}^{T}\left\Vert \bm{Y}_t-\sum_{l=1}^p\bm{A}_l\bm{Y}_{t-l}\bm{B}_l^{\top}-\bm{\Lambda}\right\Vert_F^2.
$$
Differentiating the objective function with respect to $\bm{A}_k$, $\bm{B}_k$ ($k\in\{1,...,p\}$) and $\bm{\Lambda}$ yields the following first-order conditions:
\begin{equation}\label{iteration}
    \begin{cases}
    \sum_{t=p+1}^T \left(\bm{Y}_{t}- \sum_{l=1}^p\bm{A}_l\bm{Y}_{t-l}{\bm{B}_l}^{\top}- \bm{\Lambda}\right)\bm{B}_k\bm{Y}_{t-k}^{\top}=\bm{0},~k\in\{1,...,l\},\\
    \sum_{t=p+1}^T \left(\bm{Y}_{t}^{\top}- \sum_{l=1}^p{\bm{B}_l}\bm{Y}_{t-l}^{\top}\bm{A}_l^{\top}-\bm{\Lambda}^{\top}\right)\bm{A}_k\bm{Y}_{t-k}=\bm{0},~k\in\{1,...,l\},\\
   (T-p)\bm{\Lambda}- \sum_{t=p+1}^T  \left(\bm{Y}_{t}- \sum_{l=1}^p\bm{A}_l\bm{Y}_{t-l}{\bm{B}_l}^{\top}\right)=\bm{0}.
    \end{cases}
\end{equation}

Due to the intricate product structure among the parameter matrices in (\ref{iteration}), a closed-form solution is not readily available. We therefore adopt an iterative scheme that updates one set of parameters while holding the others fixed.

\begin{enumerate}[Step 1:]
    \item \textbf{Initialization} 
    
    Set the initial values as the PROJ estimators:
    \[
    \bm{A}_l^{(0)} = \widehat{\bm{A}}_{l,PROJ},\quad \bm{B}_l^{(0)} = \widehat{\bm{B}}_{l,PROJ} \quad (l \in \{1,\dots,p\}), \quad \bm{\Lambda}^{(0)} = \widehat{\bm{\Lambda}}_{PROJ}.
    \]

    \item \textbf{Iterative updates}
    
    For each iteration $s = 0,1,2,\dots$, update the parameters sequentially using the following formulas:
    \begin{align}
        \bm{A}_l^{(s+1)} &\leftarrow \bm{M}_l^{(s)} \left( \sum_{t=p+1}^T \bm{Y}_{t-l} \bm{B}_l^{(s)^\top} \bm{B}_l^{(s)} \bm{Y}_{t-l}^\top \right)^{-1}, \nonumber \\
        \bm{B}_l^{(s+1)} &\leftarrow \bm{N}_l^{(s)} \left( \sum_{t=p+1}^T \bm{Y}_{t-l}^\top \bm{A}_l^{(s+1)^\top} \bm{A}_l^{(s+1)} \bm{Y}_{t-l} \right)^{-1}, \nonumber \\
        \bm{\Lambda}^{(s+1)} &\leftarrow \frac{1}{T-p} \sum_{t=p+1}^T \left( \bm{Y}_t - \sum_{l=1}^p \bm{A}_l^{(s+1)} \bm{Y}_{t-l} \bm{B}_l^{(s+1)^\top} \right), \nonumber
    \end{align}
    where
    \begin{align}
        \bm{M}_l^{(s)} &:= \sum_{t=p+1}^T \left( \bm{Y}_t - \sum_{k=1}^{l-1} \bm{A}_k^{(s+1)} \bm{Y}_{t-k} \bm{B}_k^{(s)^\top} - \sum_{k=l+1}^p \bm{A}_k^{(s)} \bm{Y}_{t-k} \bm{B}_k^{(s)^\top} - \bm{\Lambda}^{(s)} \right) \bm{B}_l^{(s)} \bm{Y}_{t-l}^\top, \nonumber \\
        \bm{N}_l^{(s)} &:= \sum_{t=p+1}^T \left( \bm{Y}_t^\top - \sum_{k=1}^{l-1} \bm{B}_k^{(s+1)^\top} \bm{Y}_{t-k}^\top \bm{A}_k^{(s+1)} - \sum_{k=l+1}^p \bm{B}_k^{(s)^\top} \bm{Y}_{t-k}^\top \bm{A}_k^{(s+1)} - \bm{\Lambda}^{(s)^\top} \right) \bm{A}_l^{(s+1)} \bm{Y}_{t-l}. \nonumber
    \end{align}

    \item \textbf{Convergence check}
    
    Terminate the iteration when the maximum change across all parameters falls below a predefined threshold:
    \[
    \max\left\{ \left\|\bm{A}_l^{(s+1)} - \bm{A}_l^{(s)}\right\|_F, \left\|\bm{B}_l^{(s+1)} - \bm{B}_l^{(s)}\right\|_F, \left\|\bm{\Lambda}^{(s+1)} - \bm{\Lambda}^{(s)}\right\|_F \right\} < c \times 10^{-\delta},
    \]
    for some positive constants $c$ and $\delta$.
\end{enumerate}

In this study, we set $c=1$ and $\delta=9$ without loss of generality. Upon convergence, we obtain the ICLS estimators $\widehat{\bm {A}}_{l,ICLS}$, $\widehat{\bm{B}}_{l,ICLS} ~(l\in\{1,...,p\})$ and $\widehat{\bm{\Lambda}}_{ICLS}$. The following theorem establishes the asymptotic properties of the ICLS estimators.

\begin{theorem}
Let $\{\bm{Y}_t\}$ be a MAT-NB-INAR($p$) process satisfying the conditions of Theorem \ref{thm1}. Then the ICLS estimators are strongly consistent and asymptotically normal:
$$
\sqrt{T-p}\left(\begin{array}{c}
{\rm vec}(\widehat{\bm{A}}_{1,ICLS}-\bm{A}_1) \\
{\rm vec}(\widehat{\bm{B}}_{1,ICLS}^{\top}-\bm{B}_1^{\top})\\
\vdots\\
{\rm vec}(\widehat{\bm{A}}_{p,ICLS}-\bm{A}_p) \\
{\rm vec}(\widehat{\bm{B}}_{p,ICLS}^{\top}-\bm{B}_p^{\top})\\
{\rm vec}(\widehat{\bm{V}}_{ICLS}-\bm{V})
\end{array}\right)\overset{L}{\longrightarrow}N(\bm{0},\bm{\Xi}_2),
$$
where $\bm{\Xi}_2:=\bm{Q}^{-1}E(\bm{P}_t\bm{\Sigma}_{\mathcal U}\bm{P}_t^{\top})\bm{Q}^{-1}$, with 
$\bm{Q}:=E(\bm{P}_t\bm{P}_t^{\top})+\sum_{l=1}^p\bm{\gamma}_l\bm{\gamma}_l^{\top}$, 
$$
\bm{P}_t:=\left(\bm{B}_1\bm{Y}_{t-1}^{\top}\otimes{\bm I}_m,{\bm I}_n\otimes\bm{A}_1\bm{Y}_{t-1},\cdots,\bm{B}_p\bm{Y}_{t-p}^{\top}\otimes{\bm I}_m,{\bm I}_n\otimes\bm{A}_p\bm{Y}_{t-p},{\bm I}_{mn}\right)^{\top},
$$ 
and 
$
\bm{\gamma}_l=(\bm{0}_{1\times m^2 },\bm{0}_{1\times n^2 },...,\underbrace{\bm{\alpha}_{l}^{\top},\bm{0}_{1\times n^2}}_{l\text{th position} },...,\bm{0}_{1\times m^2},\bm{0}_{1\times n^2},\bm{0}_{1\times mn})^{\top}. 
$
\end{theorem} 

\section{Simulation studies}\label{sec:6}

This section focuses on investigating the empirical performance of the proposed estimation methods. To evaluate the performance of the PROJ and ICLS methods described above, we conducted numerical simulations under three different scenarios, with sample sizes set to $T = 200$, $600$, $1000$, and performed 1000 replications.

\begin{enumerate}[\text{Scenario} A.]
  \begin{sloppypar}
  \setlength{\itemsep}{1pt}
  \setlength{\parskip}{0pt}
  \setlength{\parsep}{0pt}
  \item We consider a MAT-NB-INAR($1$) model with $(m,n)=(2,2)$, and the initial parameters in $\bm{A}$, $\bm{B}$ and $\bm{\Lambda}$ are chosen as
\[
 \widetilde{\bm{A}}= \left(\begin{array}{cc}
0.10 &0.30 \\
0.30& 0.10\\
\end{array}\right),~
\bm{B}= \left(\begin{array}{cc}
0.20 &0.40 \\
0.40& 0.20\\
\end{array}\right),~ \text{and}~
\bm{\Lambda}= \left(\begin{array}{cc}
1.00 &1.00\\
1.00&1.00\\
\end{array}\right), 
\]
and $\bm{A}=\widetilde{\bm{A}}/\Vert\widetilde{\bm{A}}\Vert_F$. 
$\bm{\mathcal{E}}_{t}$ follows a matrix-variate Poisson distribution with mean $\bm{\Lambda}$, i.e., $\bm{\mathcal{E}}_{t}\sim \text{Mpois}(\bm{\Lambda})$.

  \item We consider a MAT-NB-INAR($2$) model with $(m,n)=(2,3)$, and the initial parameters $\bm{A}_1=(a_{i,j}^{(1)})$, $\bm{A}_2=(a_{i,j}^{(2)})$, $\bm{B}_1=(b_{i,j}^{(1)})$, $\bm{B}_2=(b_{i,j}^{(2)})$ and $\bm{\Lambda}=(\lambda_{i,j})$ are chosen as
\begin{align}
 &\widetilde{\bm{A}}_1= \left(\begin{array}{cc}
0.10&0.20\\
0.40&0.50\\
\end{array}\right),~
 \widetilde{\bm{A}}_2= \left(\begin{array}{cc}
0.10&0.40\\
0.20&0.50\\
\end{array}\right),~
 \bm{B}_1= \left(\begin{array}{ccc}
0.25&0.20&0.15\\
0.20& 0.25& 0.30\\
0.15& 0.20& 0.35\\
\end{array}\right),\nonumber\\
 &\bm{B}_2= \left(\begin{array}{ccc}
0.25&0.20&0.15\\
0.20& 0.25& 0.20\\
0.15& 0.30& 0.25\\
\end{array}\right),~
\bm{\Lambda}= \left(\begin{array}{ccc}
0.50&1.50&2.00\\
2.00&1.50&0.50\\
\end{array}\right), 
\notag
\end{align}
and $\bm{A}_l=\widetilde{\bm{A}}_l/\Vert\widetilde{\bm{A}}_l\Vert_F$ ($l\in\{1,2\}$). 
$\bm{\mathcal{E}}_{t}\sim \text{MNbinom}(\bm{\Lambda},1)$.

  \item We consider a MAT-NB-INAR($2$) model with the same parameter settings as Scenario B: $(m,n)=(2,3)$; initial parameters $\bm{A}_1$, $\bm{A}_2$, $\bm{B}_1$, $\bm{B}_2$ and $\bm{\Lambda}$ are chosen identically, including the normalization $\bm{A}_l=\widetilde{\bm{A}}_l/\Vert\widetilde{\bm{A}}_l\Vert_F$ ($l\in\{1,2\}$). 
$\bm{\mathcal{E}}_t$ follows a mixed distribution, where $Z_t$ is a Bernoulli random variable with $P(Z_t = 1) = 1 - P(Z_t = 0) = p = 0.3$; $E_{1t} \sim \text{Mpois}(\bm{\Lambda})$; $E_{2t} \sim \text{MNbinom}(\bm{\Lambda}, 1)$; and $\bm{\mathcal{E}}_t = Z_t \cdot E_{1t} + (1 - Z_t) \cdot E_{2t}$.
  \end{sloppypar}
\end{enumerate}

For all the above scenarios, we rescale $\bm{A}_l=\widetilde{\bm{A}}_l/\Vert\widetilde{\bm{A}}_l\Vert_F$ before simulation to guarantee the uniqueness holds for each model. Meanwhile, all elements of $\bm{A}_l$ and $\bm{B}_l$ are chosen to satisfy $\rho({\bm{\mathcal{A}}})<1$ to guarantee the fulfillment of the stationary condition. 

\begin{figure}[h]
\centering
\includegraphics[width=5.2in]{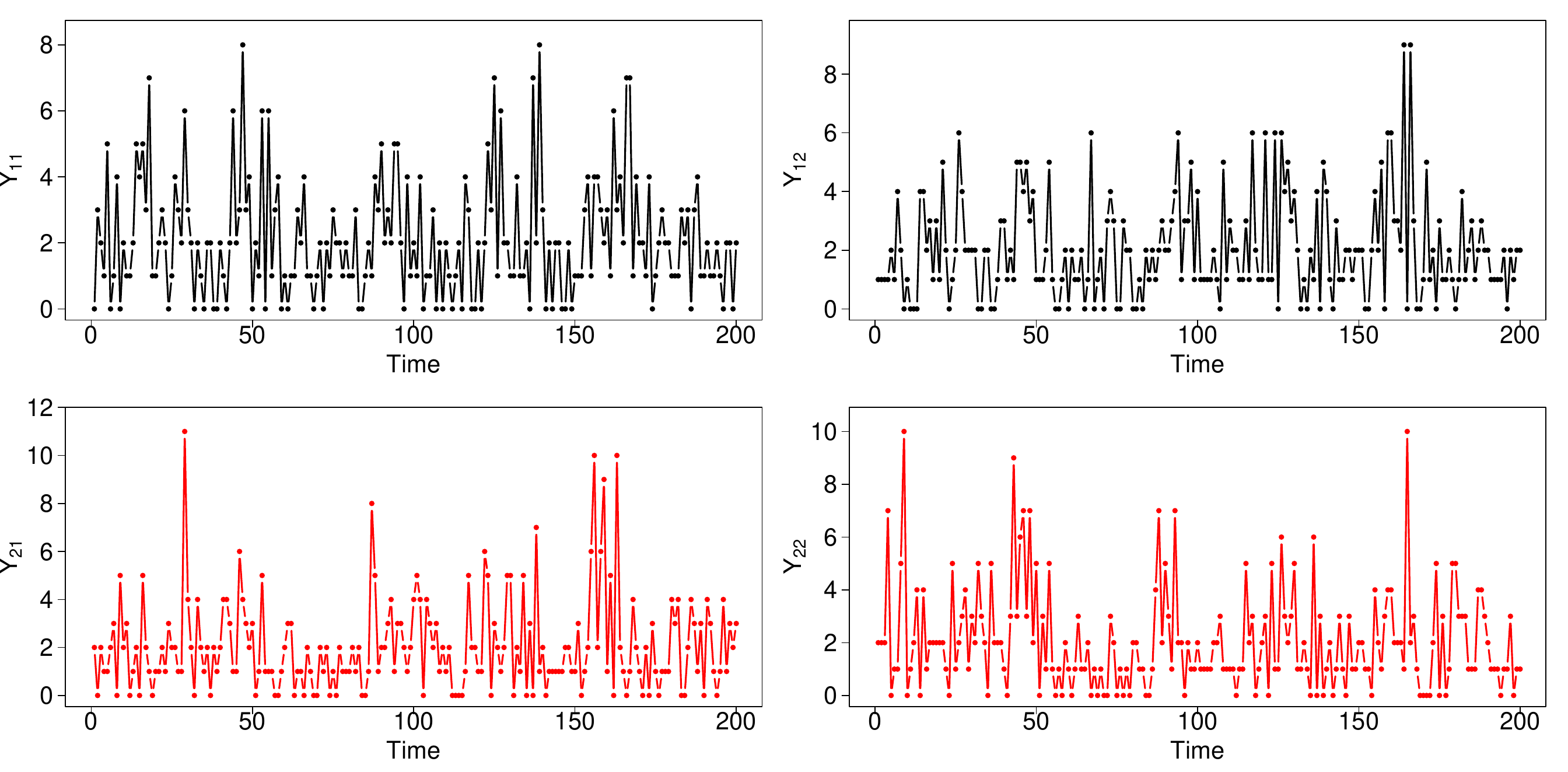}
\caption{
Time series plots of Scenario A.
}
\label{pathA}
\end{figure}
\begin{figure}[!h]
\centering
\includegraphics[width=5.6in]{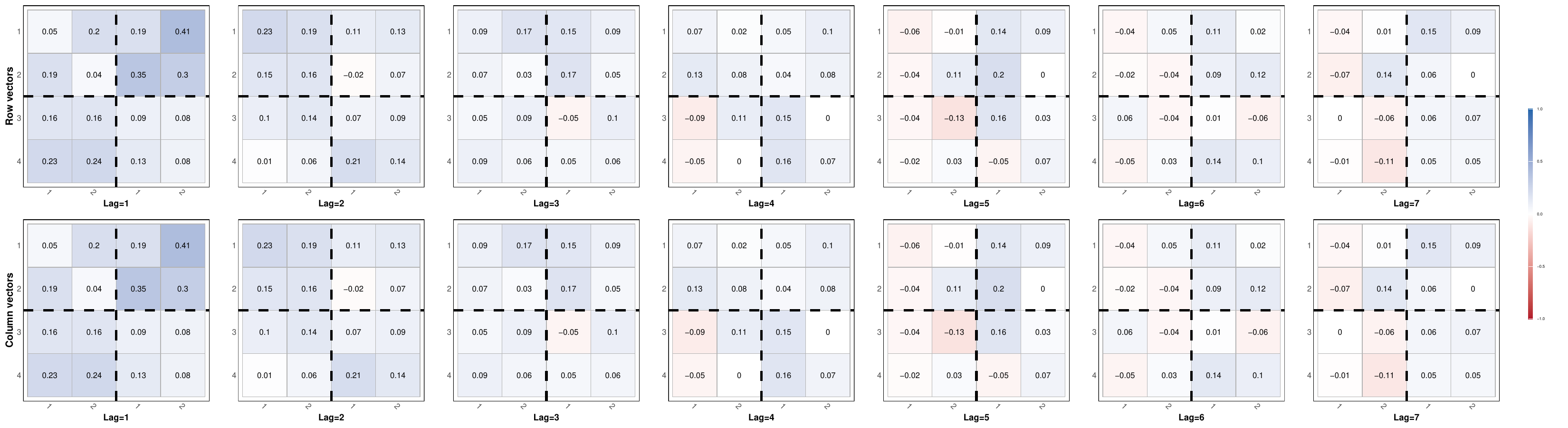}
\caption{
Cross-ACFs of row vectors and of column vectors  for Scenario A. The black bold dotted line shows the ACF matrix blocks.
}
\label{aACF}
\end{figure}

To examine the sample properties of the proposed model, we simulate a sample path under the parameter settings described above Scenario A in Fig. \ref{pathA}, displays the time series plots of $\bm{Y}_t$. As can be seen from Fig.\ref{pathA}, all four component series oscillate irregularly around their respective means throughout the entire observation window, with no apparent trend or seasonal pattern. The fluctuation range of each series remains roughly stable over time, which provides visual evidence for the stationarity of the simulated process.

To further investigate the dependence structure of the simulated matrix-variate series, we compute the cross-autocorrelation function (cross-ACF) and present the results in Fig.\ref{aACF} for lags $h=1,...,7$ when $T=200$. The time series (Fig. \ref{pathB}, \ref{pathC}) and cross-ACF plots (Fig. \ref{bACF}, \ref{cACF}) of Scenarios B and C are given in \ref{results}. Following the notation of \cite{Xu2024}, we define $\gamma_{ij,ks}(h)={\rm Cor}(y_{ij,t},y_{ks,t-h})={\rm Cov}({y}_{ij,t},{y}_{ks,t-h})/\sqrt{{\rm Var}({y}_{ij,t})}\sqrt{{\rm Var}({y}_{ks,t-h})}$, $i,k\in\{1, ... , m\}$, $j,s\in\{1,...,n\}$. The upper panel of Fig.\ref{aACF} presents the cross-ACF matrices of the row vectors of $\bm{Y}_t$, while the lower panel displays those of the column vectors. Each $4 \times 4$ matrix block corresponds to a specific lag, and is further divided into four $2 \times 2$ subblocks by the bold dashed lines, where each subblock captures the autocorrelations or cross-correlations within or between the rows (or columns) of $\bm{Y}_t$.

Several observations can be drawn from Fig.\ref{aACF}. First, at lag $h=1$, the matrix block exhibits relatively deep blue shading, with several entries attaining moderate positive values, indicating non-negligible serial dependence at short lags. As the lag $h$ increases from $1$ to $7$, the color of each matrix block progressively fades toward white and the magnitudes of the entries diminish, suggesting that the temporal dependence of the matrix series decays as the lag grows, which is consistent with the stationarity observed in Fig.\ref{pathA}. The off-diagonal subblocks in each matrix block are generally non-zero, implying that the row (or column) sequences of the proposed process exhibit not only autocorrelation within each component but also mutual cross-correlations across different positions of the matrix, which highlights the importance of adopting a matrix-variate modeling framework rather than treating the components independently. 

To evaluate the finite-sample performance of the proposed estimators, we compute the empirical bias (Bias) and standard deviation (SD) of 
the estimates across 1000 replications, together with the standard error (SE) of both the PROJ estimators and the ICLS estimators. The simulation results for Scenario A are summarized in Table \ref{t1}, while the corresponding results for Scenarios B and C are reported in Tables \ref{bt1}, \ref{bt2}, \ref{bt3} and \ref{bt4} in \ref{results}.

As shown in Table \ref{t1}, the Bias, SD, and SE of all parameter estimates decrease monotonically as the sample size $T$ increases, providing numerical evidence for the consistency of the proposed estimators. Furthermore, the SD values are in close agreement with their corresponding SE counterparts across all settings, indicating that the asymptotic approximations remain accurate even under moderate sample sizes. Regarding the comparison between the two estimators, the ICLS estimators yield smaller absolute biases than the PROJ estimators in the majority of cases, and the discrepancy between SD and SE is also more pronounced for the PROJ estimators than for the ICLS estimators. These observations collectively suggest that the ICLS estimators achieve higher estimation accuracy and better asymptotic reliability.
Parallel conclusions are drawn from the results of Scenarios B and C reported in Tables \ref{bt1}--\ref{bt4}.

\begin{landscape}
	{\captionsetup{width=1.25\textwidth}
		\begin{table}[tbp]                     
			\renewcommand\arraystretch{0.9}                  
			\setlength{\abovecaptionskip}{2pt}
			\setlength{\belowcaptionskip}{10pt}
			\caption{Simulation results for Scenario A: Bias, SE and SD}                  
			\label{t1}
			{\tabcolsep0.06in                            
				\begin{tabular}{rrrrrrrrrrrrrrr}
					\hline
					Method& $T$  & Result & $a_{1,1}$ & $a_{2,1}$ &$a_{1,2}$ & $a_{2,2}$& $b_{1,1}$ & $b_{2,1}$ &$b_{1,2}$ & $b_{2,2}$& $\lambda_{1,1}$ & $\lambda_{2,1}$ &$\lambda_{1,2}$ & $\lambda_{2,2}$\\
					\hline
					PROJ &200 &Bias&$-$0.010 &0.015 &$-$0.014 &$-$0.006 &0.001 &0.003 &0.003 &0.000 &0.018 &0.019 & 0.026 &0.021\\
					&     &SD &0.112 &0.102 &0.102 &0.111 &0.080 &0.088 &0.080 &0.088 &0.257 &0.269 &0.259 &0.262\\
					&     &SE &0.127 &0.125 &0.125 &0.127 &0.167 &0.069 &0.069 &0.167 &0.269 &0.262 &0.276 &0.276\\
					&500  &Bias&$-$0.005 &$-$0.003 &$-$0.007 &$-$0.004 &$-$0.003 &0.000 &0.004 &$-$0.002 &0.010 &0.011 &0.017 &0.009  \\
					&     & SD  & 0.072 & 0.059& 0.061& 0.070& 0.048 & 0.052 & 0.055& 0.048& 0.170& 0.166& 0.162& 0.167\\
					&     & SE  & 0.081 & 0.063& 0.066 & 0.080& 0.043 & 0.055 & 0.057& 0.043& 0.175 & 0.168 & 0.166& 0.172\\
					&1000  & Bias &$-$0.002 &0.000 & $-$0.005&$-$0.005 &0.000 &0.000 & 0.002&$-$0.001 &0.009 &0.000& 0.006&0.007  \\
					&     & SD  & 0.051& 0.043 & 0.044& 0.052& 0.034 & 0.036 & 0.037& 0.033& 0.117& 0.118 & 0.117& 0.117\\
					&     & SE  & 0.057 & 0.044& 0.046& 0.056& 0.074 & 0.070& 0.082& 0.074& 0.122 & 0.125 & 0.122& 0.122\\
					ICLS &200 &Bias&$-$0.014 &$-$0.009 &$-$0.017 &$-$0.014 &$-$0.005 &0.004 &0.004 &$-$0.001 &0.037 &0.029 &0.035 &0.017\\
					&     &SD &0.111 &0.125 &0.125 &0.111 &0.073 &0.076 &0.076 &0.073 &0.237 &0.239 &0.229 &0.230 \\
					&     &SE &0.105 &0.087 &0.088 &0.105 &0.075 &0.080 &0.080 &0.075 &0.228 &0.228 &0.228 &0.229 \\
					&500  & Bias &$-$0.004 &0.003 & $-$0.001&$-$0.006 &$-$0.005&$-$0.001 &0.003&0.001&0.013 &0.009 & 0.006&0.003\\
					&     & SD  & 0.071 & 0.085 & 0.078 & 0.071& 0.045 & 0.046 & 0.045 & 0.045 & 0.150 & 0.150 & 0.153 & 0.151\\
					&     & SE  & 0.069 & 0.058 & 0.058 & 0.069& 0.047 & 0.051& 0.051 & 0.047& 0.148 & 0.148 & 0.148 & 0.148\\
					&1000  & Bias &0.001 &$-$0.001 & 0.000&$-$0.006 &$-$0.001 &0.001 & 0.001&$-$0.001 &0.004 &0.010& 0.002&0.009  \\
					&     & SD  & 0.051 & 0.056& 0.056 & 0.049 & 0.030 & 0.032& 0.033 & 0.030 & 0.105 & 0.105& 0.105 & 0.106 \\
					&     & SE  & 0.050 & 0.041 & 0.041 & 0.049 & 0.033 & 0.036& 0.036 & 0.034 & 0.106 & 0.106& 0.106& 0.106 \\
					\hline
				\end{tabular}
			}
		\end{table}
	}
\end{landscape}

\section{Application: data}\label{sec:7}

\indent To evaluate the practical performance of the proposed matrix integer-valued autoregressive process with negative binomial thinning, we use publicly available notifiable infectious disease data from Germany, provided by the Robert Koch Institute through the SurvStat web platform:
\begin{center}
\href{http://survstat.rki.de/Content/Query/Create.aspx}{http://survstat.rki.de/Content/Query/Create.aspx}.
\end{center}
This platform supports querying case counts by disease, geographic region, and time scale, and is widely employed in infectious disease surveillance research.

Among the notifiable diseases in Germany, Campylobacter-Enteritis and Salmonellose consistently rank as the most frequently reported bacterial enteric infections. According to the German Federal Institute for Risk Assessment and the European Centre for Disease Prevention and Control, the annual number of Campylobacter-Enteritis cases typically ranges between 60000 and 70000, while Salmonellose cases lie in the range of 13000–19000. Together, these two diseases account for the vast majority of reported bacterial enteric infections. Moreover, diseases exhibit regional spatial transmission effects, and multiple pathogens often present synchronous fluctuations in seasonal co-epidemicity. Traditional univariate and conventional vector models are inadequate for characterizing their complex interactive dependency structures. For this purpose, we select two densely populated and adjacent federal states in Germany: Bavaria and Baden-Württemberg.

We selected the monthly case counts of Campylobacter-Enteritis and salmonellosis in the two federal states. The study period spans from January 2010 to December 2023, yielding a total of $T=168$ matrix-valued observations. Without loss of generality, let $y_{11,t}$ denote the monthly number of Campylobacter-Enteritis cases in Bavaria, and $y_{12,t}$ be the monthly number of Campylobacter-Enteritis cases in Baden-Württemberg. Similarly, let $y_{21,t}$  be the monthly number of Salmonellose cases in Bavaria, and $y_{22,t}$ be the monthly number of Salmonellose cases in Baden-Württemberg, where $t \in\{1,...,168\}$. Thus, at each time index $t$, the observation forms a $2\times2$ matrix-valued observation, denoted by $\bm{Y}_t$.
Among these series, the first $t\in\{1,...,134\}$ counts form the training set used to fit the model, and the last three years $t\in\{135,...,168\}$ form the testing set, which serve as the real values for $h$-step ahead out-of-sample predictions.

\begin{figure}[!h]%
\centering%
\includegraphics[width=6in]{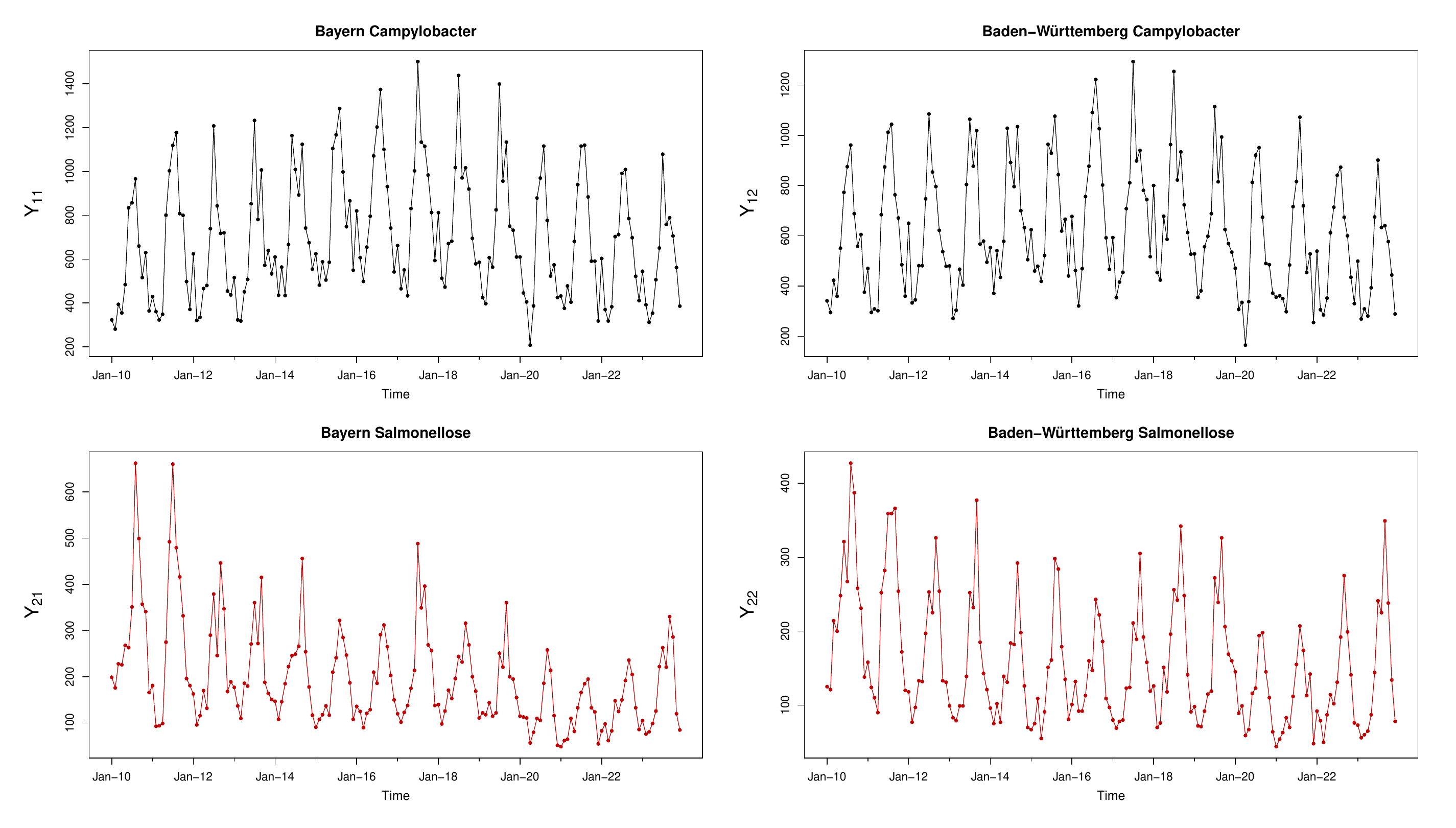}%
\caption{Time series of $\{\bm{Y}_t\}$  from January 2010 to December 2023.}
\label{trace1}
\end{figure}

Fig. \ref{trace1} displays the original monthly case counts of the two infectious diseases in Bavaria and Baden-Württemberg. None of the series display a clear trend, indicating that all analysed series are stationary. This provides an ideal data foundation for establishing a reliable predictive model in subsequent analyses.

\begin{figure}[h]
\centering
\includegraphics[width=5.2in]{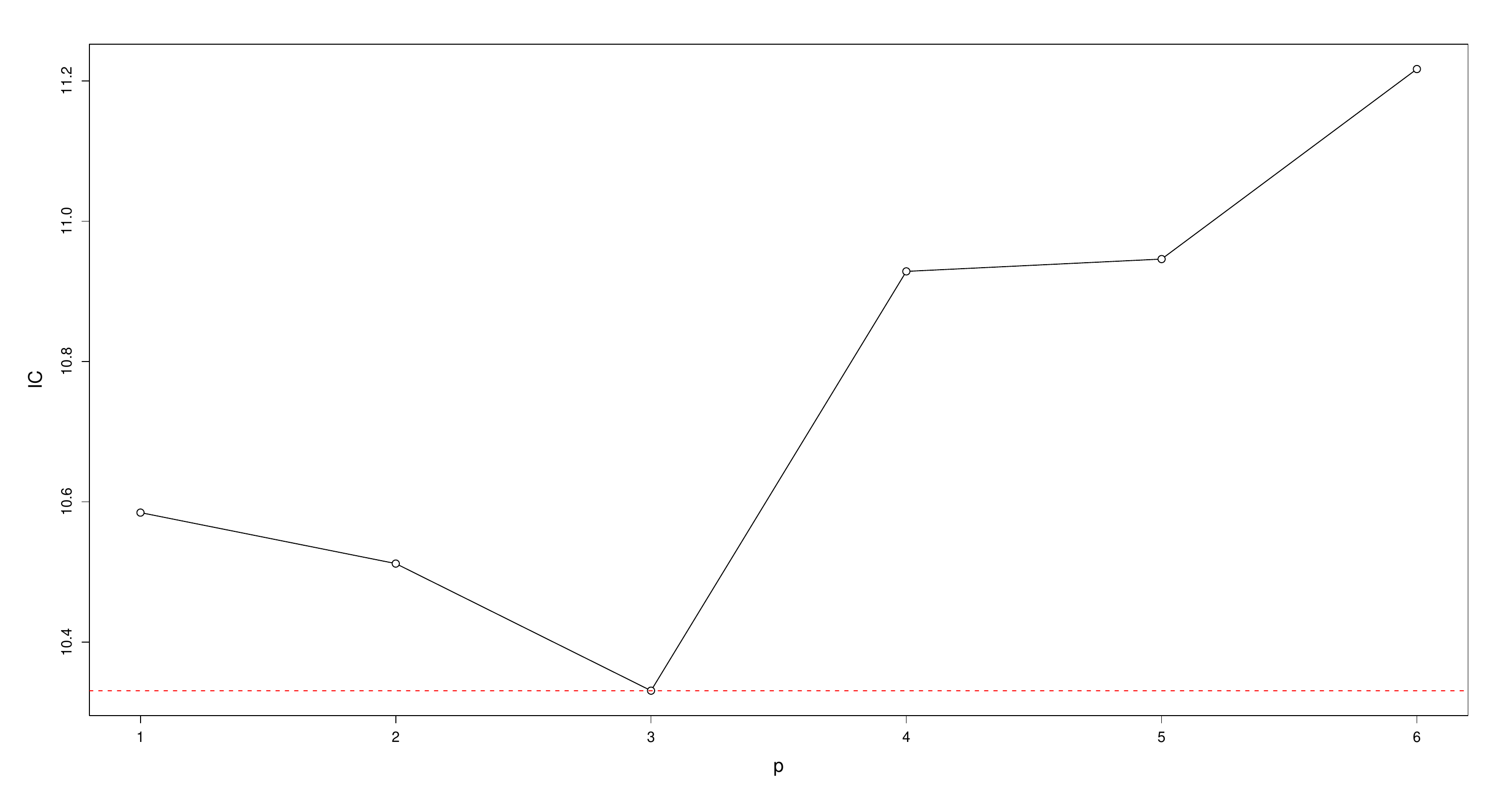}
\caption{The values of $IC_1$ with different choice $p$. The horizontal red dashed line indicates the minimum $IC_1$ value.}
\label{c_p}
\end{figure}

Then, we use the MAT-NB-INAR($p$) model and MAT-INAR($p$) to fit this data set, the order $p$ is chosen from 1 to 6. We adopt the order selection criterion proposed by \cite{Xu2024}:
\begin{equation}\label{choose_p1}
\hat{p}=\mathop{\arg\min}\limits_{1\leq \tilde p \leq \overline{p}}IC_1(\tilde p).
\end{equation}
and
\begin{equation}\label{choose_p2}
\hat{p}=\mathop{\arg\min}\limits_{1\leq \tilde p \leq \overline{p}}IC_{2}(\tilde p),
\end{equation}
To this end, we calculated $\hat{p}$ according to (\ref{choose_p1}) with $\overline{p}=6$, and draw the values of $IC_1$ in Fig. \ref{c_p}. It can be seen from the Fig. \ref{c_p} that a MAT-NB-INAR(3) model is more suitable to fit the data. For comparative reasons, we further choose the following models to fit this data set. The optimal orders of these competing models are calculated via criterion (\ref{choose_p2}):

\begin{itemize}
\item the stacked vector model, that is, MGINAR(6)  (\cite{Latour1997}) model;
\item the continuous matrix-variate autoregressive  (MAR(5))  model introduced by \cite{Chen2021};
\item the additive matrix autoregressive (Add-MAR(6)) model of \cite{Zhang2024};
\item the MAT-INAR(3) model proposed by Xu et al. (2024) using the ICLS procedure;
\item the MAT-NB-INAR($3$) model defined by (\ref{Mymodelp}) using the same ICLS estimation method.
\end{itemize}

We use the mean of residual sum of squares (MRSS) as an evaluation criterion to compare the performance of different models. First, we estimate each model using the training data, and then compute the conditional expectation $E(\bm{Y}_{t}\vert \bm{Y}_{t-1})$ as defined in equation (\ref{e1}). The MRSS of each fitting model is given as follows
$$
\text{MRSS}:=\frac{1}{T}\sum_{t=1}^T{\Big\Vert \bm{Y}_{t}-E(\bm{Y}_{t}\vert \bm{Y}_{t-1})\Big\Vert}_F.
$$
We also compare the out-of-sample forecasting performance of all models. Specifically, We predict the monthly case counts of Campylobacter-Enteritis and Salmonellose in Bavaria and Baden-Württemberg over the next three years. To evaluate the prediction accuracy, we use the mean of the out-of-sample prediction error (MSPE) as a measure, defined as follows.
$$
\text{MSPE}=\frac{1}{H}\sum_{h=1}^H{\left\Vert \widehat{\bm Y}_{t}(h)-\bm{Y}_{t+h}\right\Vert}_F,
$$
where $H=34$, $\bm{Y}_{t+h}$ denotes the true observation with $t=134$, $\widehat{\bm{Y}}_{t}(h)$ is the correspondingly predictive value obtained via the $h$-step ahead conditional expectation. All the fitting and predicted results are summarized in Table \ref{model_pk1}. We also give the number of parameters of the corresponding fitted model, expressed in $K$.

{\captionsetup{width=0.72\textwidth}
\begin{table}[tp]
\renewcommand\arraystretch{0.85}
\caption{Comparison between different models: MRSS, MSPE and $K$}
\label{model_pk1}
\vskip-0.3cm
\centering
{\tabcolsep0.3in
\begin{tabular}{cccc}
  \hline
Model               & MRSS   & MSPE &$K$     \\
  \hline
 MAT-NB-INAR(3)    &168.97   &8835.04  &28\\
 MAT-INAR(3)   &185.42  &10554.55  &28  \\
 MGINAR(6)     &125.88   &16867.78  &100 \\
Add-MAR(6)    &61.23    &16374.95  &52     \\
 MAR(5)   &208.06   &13244.32  &44  \\
\hline
\end{tabular}
}
\end{table}
}

We analyze the results in Table \ref{model_pk1} from two perspectives: in-sample fitting and out-of-sample forecasting. On the surface, the Add-MAR(6) model has the smallest MRSS among all competing models, while the MGINAR(6) model also shows a relatively good fitting performance. However, their out-of-sample forecasting performances are not satisfactory, as indicated by their relatively large MSPE values. In particular, the MGINAR(6) model involves considerably more parameters than the other models, which may lead to overfitting.

In contrast, the MAT-NB-INAR(3) model achieves the best forecasting performance among all candidate models, while maintaining a relatively small number of parameters. Compared with the MAT-INAR(3) model, the MAT-NB-INAR(3) model provides improvements in both fitting and prediction, suggesting that the negative binomial assumption is more appropriate for modeling overdispersed data. The conventional MAR(5) model performs less satisfactorily in terms of fitting, which may reflect the limitations of continuous-valued time-series models when applied to integer-valued data. Overall, the MAT-NB-INAR(3) model offers a better balance between fitting accuracy, forecasting ability, and model complexity, and is therefore the preferred model for the data under consideration.

Now, we summarize the estimated coefficient matrix results of the MAT-NB-INAR(3) model using the ICLS method in Tables \ref{model_A} -- \ref{model_V}, as well as the corresponding SE of $\widehat{\bm{A}}_{l,ICLS}$, $\widehat{\bm{B}}_{l,ICLS}~(l\in\{1,2,3\})$ and $\widehat{\bm{\Lambda}}_{ICLS}$. 

\begin{table}[htbp]
    \centering
    \renewcommand{\arraystretch}{0.9}
    \caption{Estimated left coefficient matrices $\widehat{\bm{A}}_{1,ICLS}$, $\widehat{\bm{A}}_{2,ICLS}$ and $\widehat{\bm{A}}_{3,ICLS}$ of MAT-NB-INAR(3) model.}
    \label{model_A}
    \setlength{\tabcolsep}{4pt} 
    \resizebox{\textwidth}{!}{ 
    \begin{tabular}{c cc cc cc}
        \hline
        \multirow{2}{*}{Estimation}
        & \multicolumn{2}{c}{$\widehat{\bm{A}}_{1,ICLS}$}
        & \multicolumn{2}{c}{$\widehat{\bm{A}}_{2,ICLS}$}
        & \multicolumn{2}{c}{$\widehat{\bm{A}}_{3,ICLS}$} \\
        \cline{2-3} \cline{4-5} \cline{6-7}
        & Campylobacter-Enteritis & Salmonellose
        & Campylobacter-Enteritis & Salmonellose
        & Campylobacter-Enteritis & Salmonellose \\
        \hline
        Campylobacter-Enteritis & 0.224 & 0.000 & 0.256 & 0.000 & 0.372 & 0.000 \\
        Salmonellose            & 0.000 & 0.306 & 0.000 & 0.138 & 0.000 & 0.168 \\
        \hline
    \end{tabular}
    }
\end{table}

\begin{table}[htbp]
    \centering
    \renewcommand{\arraystretch}{0.9}
    \caption{Estimated right coefficient matrices $\widehat{\bm{B}}_{1,ICLS}$, $\widehat{\bm{B}}_{2,ICLS}$ and $\widehat{\bm{B}}_{3,ICLS}$ of MAT-NB-INAR(3) model.}
    \label{model_B}
    \setlength{\tabcolsep}{4pt} 
    \resizebox{\textwidth}{!}{ 
    \begin{tabular}{c cc cc cc}
        \hline
        \multirow{2}{*}{Estimation}
        & \multicolumn{2}{c}{$\widehat{\bm{B}}_{1,ICLS}$}
        & \multicolumn{2}{c}{$\widehat{\bm{B}}_{2,ICLS}$}
        & \multicolumn{2}{c}{$\widehat{\bm{B}}_{3,ICLS}$} \\
        \cline{2-3} \cline{4-5} \cline{6-7}
        & Bayern & Baden-W\"{u}rttemberg
        & Bayern & Baden-W\"{u}rttemberg
        & Bayern & Baden-W\"{u}rttemberg \\
        \hline
        Bayern             & 0.778 & 0.096 & 0.864 & 0.190 & 0.683 & 0.474 \\
        Baden-W\"{u}rttemberg & 0.176 & 0.596 & 0.328 & 0.331 & 0.280 & 0.480 \\
        \hline
    \end{tabular}
    }
\end{table}

{\captionsetup{width=1\textwidth}
\begin{table}[t]
\renewcommand\arraystretch{0.9}
\caption{Estimated matrix $\widehat{\bm \Lambda}_{ICLS}$ of MAT-NB-INAR(3) model.}
\label{model_V}
\vskip-0.3cm
\centering
{\tabcolsep0.06in
\begin{tabular}{ccc}
 \hline
 \multirow{2}*{Estimation}  & \multicolumn{2}{c}{$\widehat{\bm \Lambda}_{ICLS}$}      \\
 \cline{2-3}
              &~~Bayern & Baden-W\"{u}rttemberg\\
  \hline
 Campylobacter-Enteritis                &~~92.260    & 224.301 \ \\
Salmonellose                    &~~85.867    & 79.171 \\
\hline
\end{tabular}
}
\end{table}
}

The estimated left and right autoregressive coefficient matrices provide several interesting findings. Table \ref{model_A} presents the estimates of $\widehat{\bm{A}}_{l,ICLS}~(l\in\{1,2\})$. Since the rows of the data matrix represent Campylobacter-Enteritis and Salmonellose, the left coefficient matrices describe the temporal relationships between the two diseases.

All three estimated left coefficient matrices are diagonal. This means that the past number of Campylobacter-Enteritis cases has little influence on the current number of Salmonellose cases, and the same is true in the opposite direction. Instead, the current number of cases of each disease mainly depends on its own past counts. Campylobacter-Enteritis has its largest coefficient at the third lag, whereas Salmonellose has its largest coefficient at the first lag. Therefore, the two diseases show different temporal patterns.

This result is reasonable because the two diseases, although both are bacterial enteric infections, do not have exactly the same sources of infection.
According to the German Federal Institute for Risk Assessment, Campylobacter infections are often associated with poultry and insufficiently cooked poultry products. Salmonella infections are more commonly associated with eggs, egg products, and raw or insufficiently cooked meat. Differences in their main sources of infection may lead to different temporal patterns and may also explain why no clear lagged relationship is found between the two diseases.

Table \ref{model_B} reports the estimated right autoregressive coefficient matrices $\widehat{\bm B}_l$ for $l\in\{1,2,3\}$. Their effects should be interpreted through $\bm B_l\circ_L\bm Y_{t-l}^{\top}$, which describes the dependence in reported disease counts between Bavaria and Baden-W\"{u}rttemberg. At the first and second lags, Bavaria has a stronger effect on Baden-W\"{u}rttemberg than Baden-W\"{u}rttemberg has on Bavaria.
At the third lag, however, the effect from Baden-W\"{u}rttemberg to Bavaria is stronger. Thus, the relationship between the two states is bidirectional, but its strength and direction change across the three lags.

From a geographical point of view, Bavaria and Baden-W\"{u}rttemberg share a long border and are closely connected by roads and railway lines.
In addition, food products are frequently transported and distributed across state borders. These movements may expose residents of the two states to similar sources of enteric infection and may help explain why changes in the case count of one state are followed by changes in the other.
The results therefore suggest a close regional relationship between Bavaria and Baden-W\"{u}rttemberg, although further mobility and food-distribution data would be needed to identify the exact transmission channels.

Table \ref{model_V} presents the estimated intercept matrix $\widehat{\bm V}$. The estimated baseline count of Campylobacter-Enteritis is higher in Baden-W\"{u}rttemberg than in Bavaria. In contrast, the baseline count of Salmonellose is slightly higher in Bavaria, although the difference between the two states is relatively small. The regional difference is therefore much more pronounced for Campylobacter-Enteritis than for Salmonellose.

These findings are also consistent with the fact that the two diseases have different major sources of infection. According to the German Federal Institute for Risk Assessment, Campylobacter infections are frequently associated with poultry and insufficiently cooked poultry products, whereas Salmonella infections are commonly related to eggs, egg products, and raw or insufficiently cooked meat. Regional differences in food production, consumption patterns, and exposure to these products may therefore be related to the different baseline levels. 

\begin{figure}[h]
\centering
\includegraphics[width=4.2in]{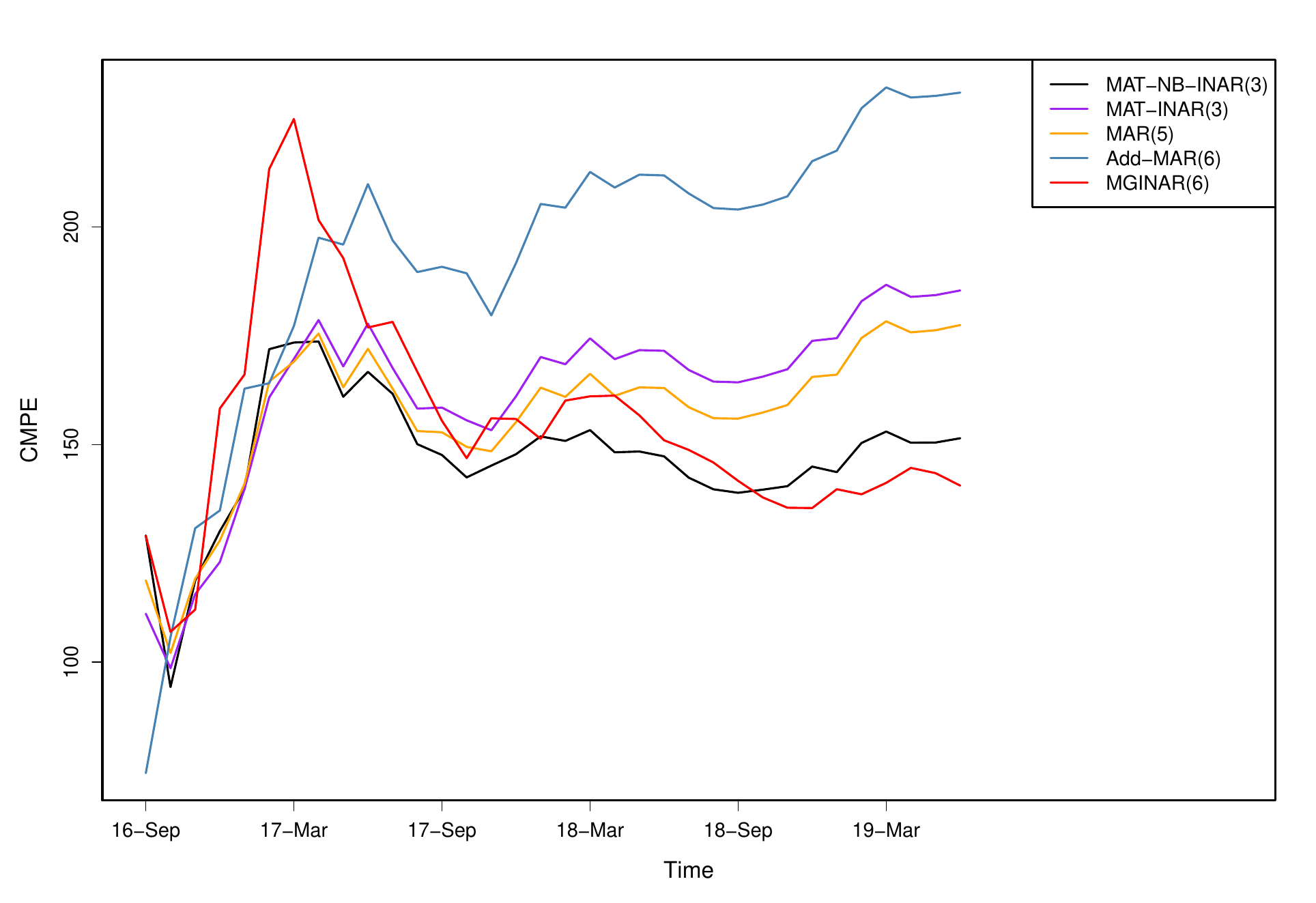}
\caption{CMPE curves of all competitive models.
}
\label{CMPE_curve}
\end{figure}

Finally, we focus on the predictive performance of the fitted  models. For the sake of comparison, we consider the cumulative loss function given by the mean-square predictive error (CMPE) as follows
$$
{\text {CMPE}}_{S}=\frac{1}{S}\sum_{h=1}^{S}{\left\Vert \widehat{\bm Y}_{t}(h)-\bm{Y}_{t+h}\right\Vert}_F, ~S\in\{1,...,36\},
$$
where $\bm{Y}_{t+h}$ denotes the true observation with $t=134$, $\widehat{\bm{Y}}_{t}(h)$ is the correspondingly predictive value obtained via the $h$-step ahead conditional expectation. The CMPE curves are presented in Fig. \ref{CMPE_curve}. As can be clearly seen in Fig. \ref{CMPE_curve}, the MAT-NB-INAR(3) model consistently achieves the lowest CMPE values across all forecast horizons, whereas the MGINAR(6) model exhibits CMPE values higher than those of all other models at most time point. Although the Add-MAR(6) model demonstrates excellent in-sample fitting performance, its predictive capability remains unsatisfactory, which further corroborates the overfitting issue arising from its excessively large number of parameters. 

It is noteworthy that the CMPE values of all models exhibit an increasing trend. This phenomenon aligns with the fundamental characteristics of infectious disease data-long-term predictions of enteric infections are subject to escalating uncertainty arising from unobserved outbreaks, seasonal fluctuations, and public health interventions. Nevertheless, the MAT-NB-INAR(3) model consistently maintains the lowest CMPE values across mostly forecast horizons, which fully demonstrates that the MAT-NB-INAR(3) model is more ideal for the prediction of matrix-variate integer-valued time series.

The fitting and prediction effects of MAT-NB-INAR(3) model for each sequence of matrix-variate data are given in Fig. \ref{fit}, the right side of the vertical black dotted line is the $h$-step ahead out-of-sample prediction effect of all comparative models with $h\in\{1,...,6\}$.

\begin{figure}[ht]
\centering
\includegraphics[width=5.6in]{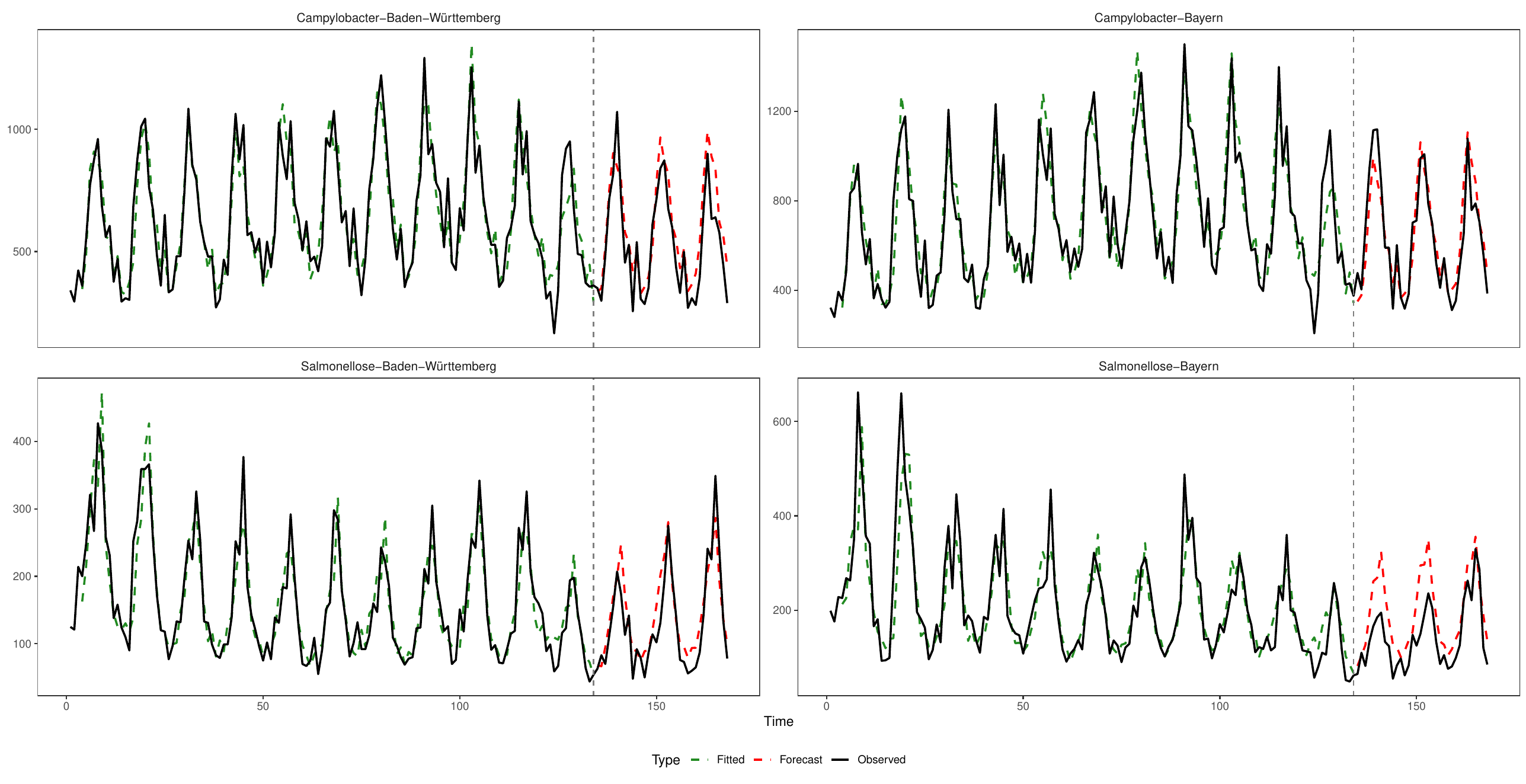}
\caption{The fitted and predicted values of the MAT-NB-INAR(3) model and the corresponding true values of the matrix-variate time series.
}
\label{fit}
\end{figure}

\section{Conclusions}\label{sec:8}

In this paper, we propose a novel matrix-valued integer-valued autoregressive process based on the negative binomial thinning operator, denoted as the MAT-NB-INAR model. The proposed model is specifically designed to address two key limitations of existing matrix autoregressive models. First, traditional continuous matrix autoregressive models are inherently unsuitable for non-negative integer-valued observations, as they fail to preserve the discrete nature of the data and cannot provide integer-valued forecasts essential for practical interpretation. Second, the recently introduced MAT-INAR model (\cite{Xu2024}), while retaining the matrix structure under binomial thinning, relies on the Poisson assumption for the innovation process. As noted in this paper, the Poisson distribution imposes equality of mean and variance, a condition frequently violated in real count data, leading to poor performance when overdispersion is present.

To overcome these drawbacks, we defined left and right matrix negative binomial thinning operators, which preserve the matrix structure while accommodating overdispersion through the geometric counting sequences induced by the thinning mechanism. We thoroughly investigated the probabilistic properties of the MAT-NB-INAR process, including stationarity, ergodicity, conditional and unconditional moments, and the cross-correlation structure between rows and columns. In terms of parameter estimation, we developed a two-stage projection method and an iterative conditional least squares procedure, and established the corresponding asymptotic normality theories. Numerical simulations provided detailed estimation results and confirmed the superiority of the ICLS method over the PROJ method.

The practical utility of the MAT-NB-INAR model is demonstrated through an application to monthly case counts of Campylobacter-Enteritis and Salmonellose in two neighbouring German federal states, Bavaria and Baden-Württemberg. The data exhibit overdispersion, cross-sectional dependence between diseases, and spatial interdependence between regions. Our model not only achieved the best performance among all competing models, but also provided interpretable coefficient matrices that reveal the temporal and spatial transmission patterns. The estimated left matrices capture the cross-disease influences, while the right matrices quantify the bidirectional regional effects that change across lags, consistent with the geographical and food-supply connections between the two states. In contrast, the MAT-INAR model based on Poisson thinning yielded larger prediction errors, confirming that the negative binomial assumption is more appropriate for such overdispersed data.

In summary, the MAT-NB-INAR model offers a coherent and powerful framework for analysing matrix-variate integer-valued time series, combining structural preservation, parameter parsimony, and the ability to handle overdispersion. This work also opens several promising directions for future research, such as extending the model to accommodate underdispersion or zero‑inflation, incorporating exogenous covariates, or developing Bayesian inference methods. We believe that our contribution will serve as a useful reference for researchers working on complex count-valued data in economics, epidemiology, criminology, and other fields where matrix-structured counts arise.

\appendix


\section{Some results of Subsection \ref{sec:7}}\label{results}
In this section, we summarize the time series and ACF plots for Scenario B in Fig. \ref{pathB} and Fig. \ref{bACF}, and for Scenario C in Fig.\ref{pathC} and Fig. \ref{cACF}. The PROJ and ICLS simulation results of Scenarios B and C are summarized in Tables \ref{bt1}--\ref{bt4}.
\begin{figure}[!h]
\centering
\includegraphics[width=5.6in]{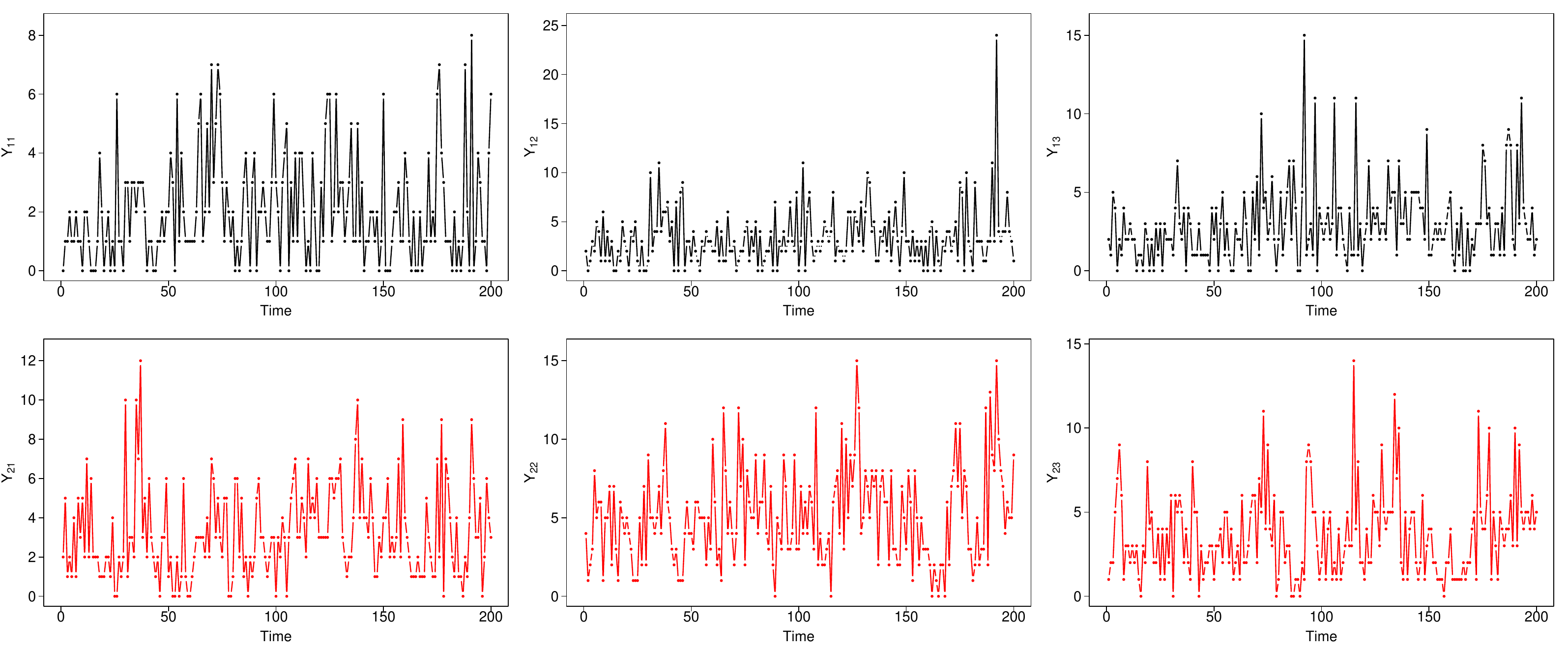}
\caption{
Time series plots for Scenario B.
}
\label{pathB}
\end{figure}
\begin{figure}[!h]
\centering
\includegraphics[width=6in]{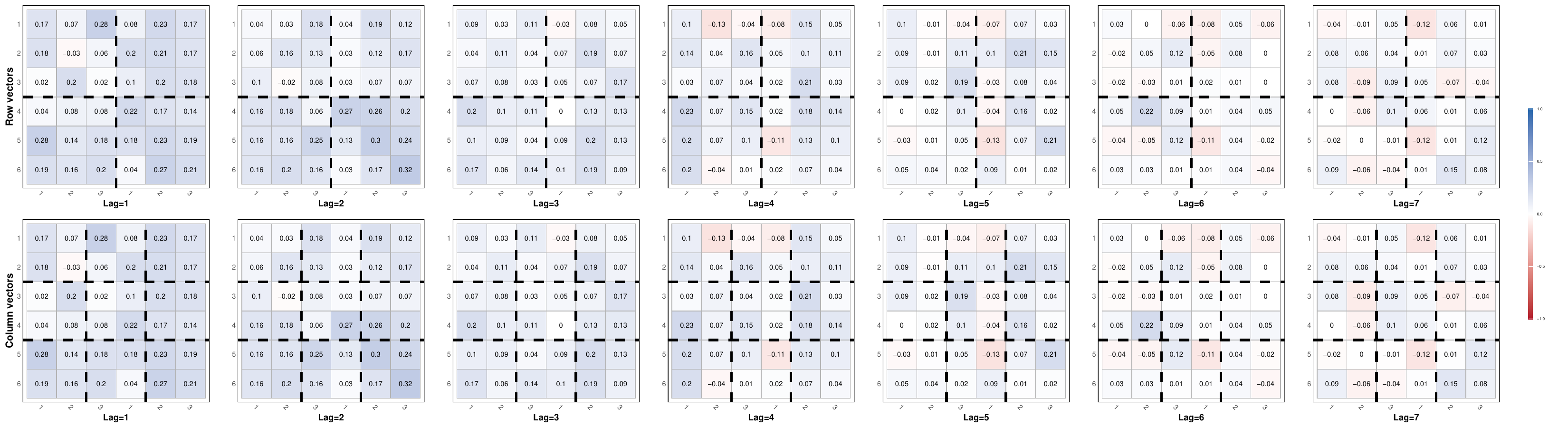}
\caption{
Cross-ACFs of row vectors (left panel) and of column vectors (right panel)  for Scenario B.
The black bold dotted line shows the ACF matrix blocks.
}
\label{bACF}
\end{figure}

\begin{figure}[h]
\centering
\includegraphics[width=6in]{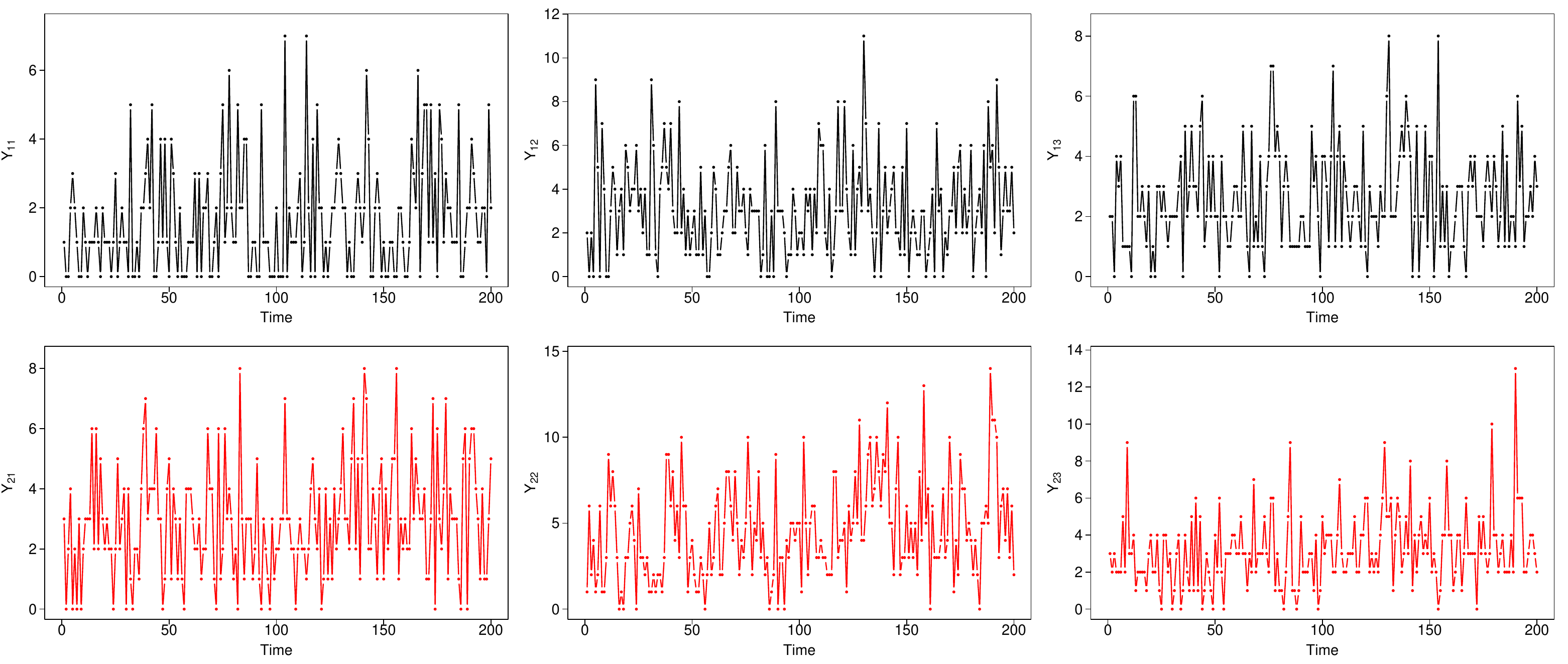}
\caption{
Time series plots for Scenario C.
}
\label{pathC}
\end{figure}

\begin{figure}[h]
\centering
\includegraphics[width=6.2in]{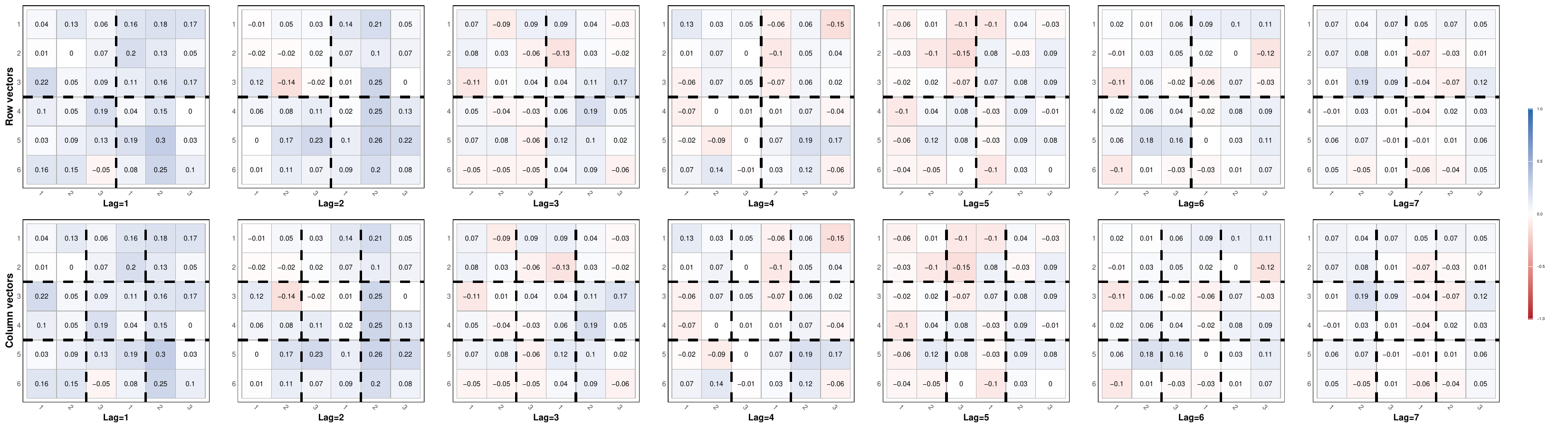}
\caption{Cross-ACFs of row vectors (left panel) and of column vectors (right panel)  for Scenario C.
The black bold dotted line shows the ACF matrix blocks.}
\label{cACF}
\end{figure}

\begin{landscape}
{\captionsetup{width=1.3\textwidth}
		\begin{table}[tbp]                     
			\renewcommand\arraystretch{0.96}                  
			\setlength{\abovecaptionskip}{2pt}
			\setlength{\belowcaptionskip}{10pt}
			\caption{PROJ simulation results for Scenario B: Bias, SE and SD}                  
			\centering                                     
			\label{bt1}
			\small
			{\tabcolsep0.05in                            
				\begin{tabular}{ccrrrrrrrrrrrrrrrr}
					\hline
					$T$  & Result & $a_{1,1}^{(1)}$ & $a_{2,1}^{(1)}$ & $a_{1,2}^{(1)}$ &$a_{2,2}^{(1)}$ & $a_{1,1}^{(2)}$ & $a_{2,1}^{(2)}$ & $a_{1,2}^{(2)}$ &$a_{2,2}^{(2)}$ & $b_{1,1}^{(1)}$ & $b_{2,1}^{(1)}$ & $b_{3,1}^{(1)}$ &$b_{1,2}^{(1)}$ & $b_{2,2}^{(1)}$ & $b_{3,2}^{(1)}$ & $b_{1,3}^{(1)}$ &$b_{2,3}^{(1)}$ \\
					\hline
					200  & Bias &0.068 &0.104 &$-$0.080&$-$0.146&0.083  &0.021&$-$0.059&$-$0.168 &0.018 &0.068 &0.049 &$-$0.006 &$-$0.024&0.095 &$-$0.036&$-$0.043\\
                    & SD   &0.186 &0.212 &0.189&0.205&0.206 &0.257 &0.180&0.235 &0.089 &0.107&0.076 &0.058 &0.081&0.077&0.063 &0.092\\
                    & SE   &0.170 &0.168 &0.130&0.130&0.183&0.165 &0.155&0.150 &0.082 &0.110 &0.088 &0.058 & 0.077&0.062 &0.073 &0.099\\
                    500  & Bias &0.046 &0.094 &$-$0.041&$-$0.074&$-$0.031 &$-$0.039&$-$0.002&$-$0.016 &0.016 &0.033 &0.025 &$-$0.009 &$-$0.023&$-$0.091 &0.024 &$-$0.065\\
                    & SD   &0.124 &0.119 &0.099&0.098&0.113&0.179&0.122&0.129&0.051 &0.061 &0.048 &0.037 &0.050&0.041 &0.038 &0.055\\
                    & SE  &0.110 &0.116&0.084&0.082&0.125 &0.124 &0.103&0.091 &0.048 &0.064 &0.050 &0.034 &0.045&0.036 &0.043 &0.058\\
                    1000  & Bias &0.037 &0.096 &$-$0.029&$-$0.061&$-$0.064 &$-$0.051 &0.012&0.021 &0.016 &0.026 &0.018 &$-$0.009 &$-$0.025&0.088 &0.020 &$-$0.066\\
                    & SD   &0.088 &0.075 &0.067&0.067&0.080 &0.115 &0.080&0.078 &0.035 &0.039&0.032&0.025 &0.033&0.028 &0.028 &0.036\\
                    & SE  &0.080 &0.084 &0.060&0.059&0.089 &0.091 &0.073&0.063 &0.033 &0.044 &0.034 &0.024 &0.032&0.025 &0.030 &0.041\\
                    \hline
                    $T$  & Result & $b_{3,3}^{(1)}$  & $b_{1,1}^{(2)}$ & $b_{2,1}^{(2)}$ & $b_{3,1}^{(2)}$ &$b_{1,2}^{(2)}$ & $b_{2,2}^{(2)}$ & $b_{3,2}^{(2)}$ & $b_{1,3}^{(2)}$ &$b_{2,3}^{(2)}$ &  $b_{3,3}^{(2)}$ & $\lambda_{1,1}$ & $\lambda_{2,1}$ & $\lambda_{1,2}$ & $\lambda_{2,2}$ & $\lambda_{1,3}$ & $\lambda_{2,3}$\\
                    \hline
                    200  & Bias&0.025 &$-$0.018 &0.060&0.058&0.005&$-$0.013 &$-$0.086&0.031 &0.097 &$-$0.037 &0.044 & 0.073&0.090 &0.163 &0.097&0.086 \\
                    & SD   &0.089 &0.089&0.112&0.081&0.062 &0.084 &0.068&0.059 &0.114 &0.077&0.333 &0.486 &0.540&0.603 &0.457 &0.459\\
                    & SE    &0.079&0.088 &0.116&0.092&0.062 &0.082 &0.065&0.074 &0.098 &0.078 &0.372 &0.522 &0.586&0.663 &0.496 &0.514\\
                    500  & Bias &0.026 &$-$0.028 &0.005&0.011&0.009 &0.003 &$-$0.079&0.001&0.080 &$-$0.050 &0.019 &0.046&0.038&0.073 &0.054 &0.060\\
                    & SD   &0.050 &0.057 &0.069&0.044&0.039 &0.050 &0.041&0.039 &0.072 &0.054 &0.214&0.308 &0.341&0.387 &0.289 &0.308\\
                    & SE    &0.046 &0.050 &0.065&0.050&0.037 &0.048 &0.037&0.044&0.058 &0.045 &0.220 &0.314 &0.347&0.400 &0.292 &0.306\\
                    1000  & Bias  &0.025 &$-$0.020 &$-$0.011&$-$0.008&0.012 &0.014 &$-$0.077 &$-$0.012 &0.076 &$-$0.045&0.004&0.013 &0.026&0.036 &0.053 &0.043\\
                    & SD    &0.033 &0.039 &0.049&0.029&0.025 &0.034 &0.028&0.027 &0.049 &0.039 &0.155 &0.225 &0.239&0.277 &0.212 &0.219\\
                    & SE    &0.032 &0.034 &0.045&0.034&0.026 &0.034 &0.026&0.031 &0.040 &0.031 &0.153&0.220 &0.243&0.281 &0.203 &0.214\\
					\hline
				\end{tabular}
			}
		\end{table}
	}
	{\captionsetup{width=1.3\textwidth}
		\begin{table}[tbp]                     
			\renewcommand\arraystretch{0.96}                  
			\setlength{\abovecaptionskip}{2pt}
			\setlength{\belowcaptionskip}{10pt}
			\caption{ICLS simulation results for Scenario B: Bias, SE and SD}                  
			\centering                                     
			\label{bt2}
			\small
			{\tabcolsep0.05in                            
				\begin{tabular}{ccrrrrrrrrrrrrrrrr}
					\hline
					$T$  & Result & $a_{1,1}^{(1)}$ & $a_{2,1}^{(1)}$ & $a_{1,2}^{(1)}$ &$a_{2,2}^{(1)}$ & $a_{1,1}^{(2)}$ & $a_{2,1}^{(2)}$ & $a_{1,2}^{(2)}$ &$a_{2,2}^{(2)}$ & $b_{1,1}^{(1)}$ & $b_{2,1}^{(1)}$ & $b_{3,1}^{(1)}$ &$b_{1,2}^{(1)}$ & $b_{2,2}^{(1)}$ & $b_{3,2}^{(1)}$ & $b_{1,3}^{(1)}$ &$b_{2,3}^{(1)}$ \\
					\hline
					200  & Bias &$-$0.038 &$-$0.035 &$-$0.084&$-$0.100&$-$0.047  &$-$0.039&$-$0.032&$-$0.083 &0.012&0.013 &0.009&0.016 &0.016&0.129&0.010 &$-$0.077\\
					& SD   &0.183 &0.218 &0.155&0.177&0.207&0.253 &0.170&0.204 &0.089 &0.125&0.095 &0.064 &0.086&0.072&0.086 &0.128\\
					& SE   &0.136 &0.149 &0.105&0.099&0.152 &0.151 &0.126&0.120 &0.078 &0.104&0.082 &0.060& 0.079&0.064 &0.076 &0.102\\
					500  & Bias &$-$0.021 &$-$0.026 &$-$0.046&$-$0.046&$-$0.006 &0.016&0.009&0.010 &0.010 &0.010 &0.004 &0.006 &0.009&0.115 &0.003 &$-$0.094\\
					& SD   &0.115 &0.136 &0.095&0.113&0.126&0.161&0.116&0.138&0.050 &0.062 &0.051 &0.036 &0.047&0.038 &0.047 &0.061\\
					& SE  &0.101 &0.112 &0.077&0.069&0.118 &0.120 &0.098&0.092 &0.047 &0.062 &0.049&0.035 &0.047&0.038 &0.046 &0.062\\
					1000  & Bias &$-$0.011 &$-$0.021 &$-$0.030&$-$0.039&$-$0.004 &0.023 &0.014&0.022 &0.007 &0.006&0.002 &0.005&0.007&0.112 &0.002 &$-$0.095\\
					& SD   &0.080 &0.098&0.066&0.080&0.091 &0.117&0.081&0.100 &0.034 &0.043&0.035&0.023 &0.032&0.026 &0.034 &0.042\\
					& SE  &0.075 &0.084 &0.056&0.052 &0.086 &0.088&0.072 &0.067 &0.033&0.043 &0.034 &0.025 &0.033&0.027 &0.032 &0.043\\
					\hline
					$T$  & Result & $b_{3,3}^{(1)}$  & $b_{1,1}^{(2)}$ & $b_{2,1}^{(2)}$ & $b_{3,1}^{(2)}$ &$b_{1,2}^{(2)}$ & $b_{2,2}^{(2)}$ & $b_{3,2}^{(2)}$ & $b_{1,3}^{(2)}$ &$b_{2,3}^{(2)}$ &  $b_{3,3}^{(2)}$ & $\lambda_{1,1}$ & $\lambda_{2,1}$ & $\lambda_{1,2}$ & $\lambda_{2,2}$ & $\lambda_{1,3}$ & $\lambda_{2,3}$\\
					\hline
					200  & Bias&0.009 &$-$0.005&0.009&0.010&0.007 &$-$0.012 &$-$0.089&0.005 &0.116 &$-$0.006 &0.095 & 0.093&0.125 &0.119&0.102&0.105\\
					& SD   &0.090 &0.094&0.124&0.106&0.070 &0.091 &0.074&0.089 &0.117 &0.091&0.270 &0.422 &0.426&0.537 &0.372 &0.429\\
					& SE    &0.081&0.081 &0.106&0.082&0.061 &0.080 &0.063&0.073 &0.098 &0.076 &0.280 &0.408 &0.400&0.536 &0.340 &0.413\\
					500  & Bias &0.009 &$-$0.008 &$-$0.002&0.003&$-$0.003 &$-$0.011 &$-$0.103&$-$0.000 &0.095 &$-$0.008 &0.041&0.033&0.044&0.044&0.050 &0.033\\
					& SD   &0.048 &0.048 &0.062&0.047&0.036 &0.045 &0.038&0.044 &0.057 &0.046 &0.180&0.266&0.255&0.345 &0.229&0.270\\
					& SE    &0.049 &0.047 &0.061&0.047&0.036 &0.046 &0.036&0.043&0.057 &0.044 &0.179 &0.263 &0.258&0.348 &0.220 &0.269\\
					1000  & Bias  &0.008 &$-$0.008&$-$0.006&$-$0.002&$-$0.003 &$-$0.009 &$-$0.104 &$-$0.002&0.093&$-$0.007&0.021&0.028&0.024&0.032 &0.018&0.024\\
					& SD    &0.034 &0.034 &0.042&0.032&0.026 &0.030 &0.025&0.030 &0.038&0.031 &0.130 &0.194 &0.176&0.245 &0.155 &0.190\\
					& SE    &0.034 &0.033 &0.042&0.033&0.025 &0.033 &0.026&0.030 &0.040 &0.031 &0.128&0.188 &0.183&0.248&0.156 &0.192\\
					\hline
				\end{tabular}
			}
		\end{table}
	}

\end{landscape}


\begin{landscape}
	{\captionsetup{width=1.3\textwidth}
		\begin{table}[tbp]                     
			\renewcommand\arraystretch{0.96}                  
			\setlength{\abovecaptionskip}{2pt}
			\setlength{\belowcaptionskip}{10pt}
			\caption{PROJ simulation results for Scenario C: Bias, SE and SD}                  
			\centering                                     
			\label{bt3}
			\small
			{\tabcolsep0.05in                            
				\begin{tabular}{ccrrrrrrrrrrrrrrrr}
					\hline
					$T$  & Result & $a_{1,1}^{(1)}$ & $a_{2,1}^{(1)}$ & $a_{1,2}^{(1)}$ &$a_{2,2}^{(1)}$ & $a_{1,1}^{(2)}$ & $a_{2,1}^{(2)}$ & $a_{1,2}^{(2)}$ &$a_{2,2}^{(2)}$ & $b_{1,1}^{(1)}$ & $b_{2,1}^{(1)}$ & $b_{3,1}^{(1)}$ &$b_{1,2}^{(1)}$ & $b_{2,2}^{(1)}$ & $b_{3,2}^{(1)}$ & $b_{1,3}^{(1)}$ &$b_{2,3}^{(1)}$ \\
					\hline
					200  & Bias &0.071 &0.105 &$-$0.084&$-$0.147&0.084  &0.016 &$-$0.063&$-$0.170 &0.023 &0.061 &0.054 &0.001 &$-$0.018&0.096 &0.032 &$-$0.056\\
					& SD   &0.192 &0.218 &0.188&0.199&0.210 &0.264 &0.184&0.242 &0.090 &0.101&0.077 &0.062 &0.080&0.086&0.060 &0.091\\
					& SE   &0.168 &0.168 &0.128&0.131&0.179 &0.164 &0.151&0.148 &0.083 &0.108&0.090 &0.060 & 0.078&0.064 &0.074 &0.096\\
					500  & Bias &0.040 &0.103 &$-$0.038&$-$0.079&$-$0.038 &$-$0.043&$-$0.000&$-$0.010 &0.016 &0.031 &0.022 &$-$0.007 &$-$0.025&$-$0.090 &0.022 &$-$0.067\\
					& SD   &0.120 &0.117&0.101&0.099&0.111&0.173&0.116&0.126&0.052 &0.060&0.044 &0.036 &0.051&0.043 &0.040&0.051\\
					& SE  &0.110 &0.115 &0.083&0.083&0.122&0.123 &0.100&0.089 &0.048 &0.062 &0.051 &0.035 &0.046&0.037 &0.043 &0.057\\
					1000  & Bias &0.041 &0.097 &$-$0.026&$-$0.064&$-$0.073 &$-$0.062 &0.010&0.033 &0.015 &0.024 &0.019 &$-$0.009&$-$0.024&0.088 &0.018 &$-$0.071\\
					& SD   &0.085 &0.071 &0.068&0.062&0.075 &0.113 &0.077&0.072&0.034 &0.040&0.033&0.025 &0.035&0.028 &0.027 &0.035\\
					& SE  &0.079 &0.084 &0.059&0.059&0.088 &0.092 &0.071&0.061 &0.033 &0.043 &0.035 &0.024 &0.032&0.025 &0.030 &0.039\\
					\hline
					$T$  & Result & $b_{3,3}^{(1)}$  & $b_{1,1}^{(2)}$ & $b_{2,1}^{(2)}$ & $b_{3,1}^{(2)}$ &$b_{1,2}^{(2)}$ & $b_{2,2}^{(2)}$ & $b_{3,2}^{(2)}$ & $b_{1,3}^{(2)}$ &$b_{2,3}^{(2)}$ &  $b_{3,3}^{(2)}$ & $\lambda_{1,1}$ & $\lambda_{2,1}$ & $\lambda_{1,2}$ & $\lambda_{2,2}$ & $\lambda_{1,3}$ & $\lambda_{2,3}$\\
					\hline
					200  & Bias&0.026 &$-$0.015 &0.054&0.057&0.003 &$-$0.015 &$-$0.079&0.031 &0.097 &$-$0.037 &0.025 & 0.085&0.093 &0.166 &0.091&0.104 \\
					& SD   &0.090 &0.092&0.110&0.082&0.062 &0.087 &0.068&0.062 &0.109 &0.078 &0.337 &0.483 &0.515&0.600 &0.457 &0.486\\
					& SE    &0.080&0.088 &0.114&0.093&0.064 &0.082 &0.068&0.075 &0.096 &0.079 &0.378 &0.525 &0.570&0.652 &0.493 &0.526\\
					500  & Bias &0.024 &$-$0.025 &0.000&0.009&0.012 &0.005 &$-$0.080&0.001 &0.083 &$-$0.042 &0.013 &0.018&0.053&0.081 &0.051 &0.060\\
					& SD   &0.049 &0.055 &0.066&0.042&0.040 &0.050 &0.041&0.036 &0.069 &0.055 &0.218&0.303 &0.322&0.371&0.275 &0.289\\
					& SE    &0.046 &0.049 &0.064&0.051&0.037 &0.048 &0.038&0.044&0.057 &0.045 &0.223 &0.313&0.336&0.395 &0.290 &0.312\\
					1000  & Bias  &0.022 &$-$0.020 &$-$0.010&$-$0.006&0.011 &0.010 &$-$0.078 &$-$0.010 &0.078 &$-$0.039&0.007&0.019 &0.043&0.038 &0.041 &0.033\\
					& SD    &0.033 &0.037 &0.048&0.031&0.027 &0.035 &0.028&0.028 &0.047 &0.038 &0.158 &0.226 &0.237&0.268 &0.199 &0.208\\
					& SE    &0.032 &0.034 &0.044&0.035&0.026 &0.034 &0.027&0.030 &0.039 &0.031 &0.155&0.219 &0.235&0.277 &0.202 &0.218\\
					\hline
				\end{tabular}
			}
		\end{table}
	}
	{\captionsetup{width=1.3\textwidth}
		\begin{table}[tbp]                     
			\renewcommand\arraystretch{0.96}                  
			\setlength{\abovecaptionskip}{2pt}
			\setlength{\belowcaptionskip}{10pt}
			\caption{ICLS simulation results for Scenario C: Bias, SE and SD}                  
			\centering                                     
			\label{bt4}
			\small
			{\tabcolsep0.05in                            
				\begin{tabular}{ccrrrrrrrrrrrrrrrr}
					\hline
					$T$  & Result & $a_{1,1}^{(1)}$ & $a_{2,1}^{(1)}$ & $a_{1,2}^{(1)}$ &$a_{2,2}^{(1)}$ & $a_{1,1}^{(2)}$ & $a_{2,1}^{(2)}$ & $a_{1,2}^{(2)}$ &$a_{2,2}^{(2)}$ & $b_{1,1}^{(1)}$ & $b_{2,1}^{(1)}$ & $b_{3,1}^{(1)}$ &$b_{1,2}^{(1)}$ & $b_{2,2}^{(1)}$ & $b_{3,2}^{(1)}$ & $b_{1,3}^{(1)}$ &$b_{2,3}^{(1)}$ \\
					\hline
					200  & Bias &$-$0.034 &$-$0.040 &$-$0.086&$-$0.108&$-$0.046  &$-$0.037&$-$0.036&$-$0.078 &0.012 &0.018 &0.006 &0.014 &0.018&0.127 &0.010 &$-$0.083\\
					& SD   &0.178 &0.236 &0.152&0.178&0.192 &0.259 &0.174&0.210 &0.087 &0.120&0.093 &0.064 &0.089&0.071&0.085 &0.113\\
					& SE   &0.138 &0.151 &0.105&0.100&0.153 &0.153 &0.125&0.120 &0.079 &0.103 &0.084 &0.060 & 0.078&0.066 &0.075 &0.099\\
					500  & Bias &$-$0.016 &$-$0.025 &$-$0.038&$-$0.048&$-$0.012 &0.010&0.007&0.013 &0.006 &0.006 &0.003 &0.008 &0.010&0.113 &0.003 &$-$0.096\\
					& SD   &0.120 &0.140 &0.099&0.111&0.131&0.165&0.110&0.137&0.049 &0.061 &0.050 &0.036 &0.045&0.036 &0.046 &0.060\\
					& SE  &0.102 &0.116 &0.077&0.071&0.116&0.119 &0.095&0.090 &0.047 &0.061 &0.049&0.036 &0.047&0.039 &0.045 &0.060\\
					1000  & Bias &$-$0.009 &$-$0.013 &$-$0.029&$-$0.033&$-$0.002&0.022&0.014&0.018 &0.006 &0.005 &0.001&0.005&0.007&0.110 &0.002 &$-$0.098\\
					& SD   &0.077 &0.100 &0.068&0.082&0.091 &0.116&0.080&0.097&0.033 &0.043&0.035&0.025 &0.031&0.026 &0.032 &0.040\\
					& SE  &0.076 &0.086 &0.057&0.052 &0.086 &0.088&0.071 &0.067&0.033&0.042 &0.034 &0.025 &0.033&0.028 &0.031 &0.042\\
					\hline
					$T$  & Result & $b_{3,3}^{(1)}$  & $b_{1,1}^{(2)}$ & $b_{2,1}^{(2)}$ & $b_{3,1}^{(2)}$ &$b_{1,2}^{(2)}$ & $b_{2,2}^{(2)}$ & $b_{3,2}^{(2)}$ & $b_{1,3}^{(2)}$ &$b_{2,3}^{(2)}$ &  $b_{3,3}^{(2)}$ & $\lambda_{1,1}$ & $\lambda_{2,1}$ & $\lambda_{1,2}$ & $\lambda_{2,2}$ & $\lambda_{1,3}$ & $\lambda_{2,3}$\\
					\hline
					200  & Bias&0.016 &$-$0.001 &0.001&0.004&0.004 &$-$0.006 &$-$0.090&0.005 &0.116 &0.000 &0.108 & 0.108&0.109 &0.110 &0.116&0.098 \\
					& SD   &0.090 &0.091&0.125&0.092&0.068 &0.092 &0.073&0.082 &0.112 &0.088 &0.276 &0.411 &0.405&0.526 &0.355 &0.425\\
					& SE    &0.081&0.080 &0.103&0.083&0.062 &0.080 &0.065&0.073 &0.095 &0.077 &0.278 &0.406 &0.391&0.526 &0.337 &0.422\\
					500  & Bias &0.007 &$-$0.009 &$-$0.003&$-$0.001&$-$0.002 &$-$0.008 &$-$0.099&$-$0.002 &0.095 &$-$0.008 &0.035 &0.035&0.046&0.052 &0.028 &0.029\\
					& SD   &0.051 &0.049 &0.061&0.049&0.036 &0.044 &0.038&0.045 &0.057 &0.046 &0.181&0.270 &0.253&0.343 &0.222 &0.274\\
					& SE    &0.049 &0.047 &0.060&0.048&0.036 &0.047 &0.038&0.043&0.056&0.045 &0.180 &0.262 &0.252&0.343 &0.218 &0.275\\
					1000  & Bias  &0.007 &$-$0.009 &$-$0.004&0.000&$-$0.004 &$-$0.008 &$-$0.104 &$-$0.002&0.093 &$-$0.007&0.022&0.024 &0.028&0.027&0.022 &0.015\\
					& SD    &0.031 &0.034 &0.041&0.032&0.026 &0.030 &0.026&0.031 &0.038&0.032 &0.129 &0.188 &0.184&0.237 &0.162 &0.199\\
					& SE    &0.034 &0.033 &0.042&0.034&0.026 &0.033 &0.027&0.030 &0.039 &0.032 &0.128&0.187 &0.179&0.245&0.155 &0.196\\
					\hline
				\end{tabular}
			}
		\end{table}
	}
	
\end{landscape}

\bibliography{mybibfile}
\end{document}